\documentclass[10pt]{amsart}
\RequirePackage{graphicx}
\usepackage[top=3cm,bottom=2cm,left=3cm,right=3cm,marginparwidth=1.75cm]{geometry}
\usepackage{mathrsfs}   
\usepackage[english]{babel}
\usepackage[T1]{fontenc}
\usepackage{dsfont}
\usepackage{xcolor, import}
\usepackage{dsfont}

\usepackage{fancyvrb}
\usepackage[colorlinks=true, citecolor=red, linkcolor=blue]{hyperref}
\usepackage{xparse}
\usepackage[titletoc,title]{appendix}
\usepackage{amsfonts,mathtools}
\usepackage{amsmath,amssymb,amsthm}
\usepackage{bbm,graphicx,xargs}
\usepackage[square,compress,comma, numbers,sort]{natbib}
\usepackage{mathbbol}
\usepackage[shortlabels]{enumitem}
\usepackage{cleveref}

\numberwithin{equation}{section}
 \newtheorem{theo}{Theorem}[section]
 \newtheorem{sat}[theo]{Proposition}
 \newtheorem{de}[theo]{Definition}
 \newtheorem{lem}[theo]{Lemma}
 
 \newtheorem{example}[theo]{Example}
 \newtheorem{korr}[theo]{Corollary}
 \newtheorem{remark}[theo]{Remark}
 \newcommand{\BEL}{\begin{lem}}
	\newcommand{\EEL}{\end{lem}}

\newcommand{\BRM}{\begin{remark}}
	\newcommand{\ERM}{\end{remark}}

\newcommand{\BD}{\begin{de}}
	\newcommand{\ED}{\end{de}}
\newcommand{\BS}{\begin{sat}}
	\newcommand{\ES}{\end{sat}}
\newcommand{\BT}{\begin{theo}}
	\newcommand{\ET}{\end{theo}}
\newcommand{\BK}{\begin{korr}}
	\newcommand{\EK}{\end{korr}}

\def\Ife{\mathcal{I}_{\EI,fe}}  
\def\Ifasup{ \mathcal{I}_{\EI,argmax}}   

\newcommand{\R}{\mathbb{R}}
\newcommand{\inr}{\in \R}
\newcommand{\QED}{\hfill $\Box$}

\DeclarePairedDelimiter{\parens}()
\def\Pspace{ \parens{\Omega,\mathscr{F},\mathbb{P}}}

\def\borel#1{\mathscr{B}(#1)}

\def\AA{\mathscr{D}}

\def\TTd{\EIE }

\newcommand{\abs}[1]{\left\lvert #1 \right\rvert}

\DeclarePairedDelimiterXPP\pk[1]{\mathbb{P}}\{ \}{}{ #1}
\DeclarePairedDelimiterXPP\pkh[1]{\widehat{\mathbb{P}}}\{ \}{}{ #1}
\DeclarePairedDelimiterXPP\pkT[1]{\widetilde{\mathbb{P}}}\{ \}{}{ #1}
\DeclarePairedDelimiterXPP\E[1]{\mathbb{E}}\{ \}{}{	#1}
\DeclarePairedDelimiterXPP\EHh[1]{\widehat{\mathbb{E}}}\{ \}{}{	#1}

\def\kK{\mathscr{C}_{\kappa} [\vk Z]}
\def\kkA{\mathcal{C}_{\kappa}[Z_\kappa^{1/\alpha}]}
\def\kk{\mathcal{C}_{\kappa} [\vk Z ]}
\def\kkQ{\mathcal{C}_{\kappa,N}[  \vk Q ]}
\def\kkQm{\mathcal{C}_{\kappa, N}[  \vk Q^{(m)} ]}
\def\kkL{\mathcal{C}_{\kappa,N}[  \vk Q ]}

\def\AA{\mathscr{D}}

\def\II{I}
\def\bobo{\mathbb{B}}

\def\SZq{ \mathcal{S}_{\EI} (\vk Z)}
\def\STq{ \mathcal{S}_{\EI} (\vk \Theta)}
\def\JJq{ \mathfrak{B}_{\EI,\tau} (\vk Y)}

\def\JJxi{ \mathcal{B,\tau} (\vk Y  )}

\def\SQ{ \mathcal{S} (\vk Q)}

\def\IF{\infty}

\newcommand{\BQN}{\begin{eqnarray}}
	\newcommand{\EQN}{\end{eqnarray}}
\newcommand{\BQNY}{\begin{eqnarray*}}
	\newcommand{\EQNY}{\end{eqnarray*}}
\def\ldot{, \ldots,}

\newcommand{\limit}[1]{\lim_{#1 \to   \infty}}

\def\bqn#1{ \begin{eqnarray} #1 \end{eqnarray}}
\def\bqny#1{ \begin{eqnarray*} #1 \end{eqnarray*}}

\newcommand{\kb}[1]{\boldsymbol{#1}}
\newcommand{\vk}[1]{\kb{#1}}

\newcommand{\COM}[1]{}

\definecolor{c20}{rgb}{0.,0.7,0.}
\definecolor{c30}{rgb}{0.,0.,1.}
\definecolor{c40}{rgb}{1,0.1,0.7}
\definecolor{c50}{rgb}{1,0,0}
\definecolor{c60}{rgb}{1,0.9,0.1}
\definecolor{c70}{rgb}{0.50,1.00,0.00}

\def\cED#1{{\textcolor{c50}{#1}}}
\def\cEE#1{{#1}}

\def\cF#1{{\textcolor{c50}{#1}}}
\def\cF#1{{#1}}
\def\cEE#1{{#1}}

\def\cED#1{#1}

\newcommand\Dset{\mathfrak{D}}
\newcommand\DDset{\Dset_C}

\def\intDD{\int_{\EI}}

\def\Hh{ \mathfrak{H}}

\def\clasP{\mathfrak{W}_{\kappa}}

\def\MPY{ \mathcal{M}_\EI(\vk Y)}

\def\SY{\mathcal{S}(\vk Y)}

\DeclarePairedDelimiterXPP\ETT[1]{\widetilde{\mathbb{E}}}\{ \}{}{	#1}

\DeclarePairedDelimiterXPP\ind[1]{\mathbb{I}}( ){}{	#1}

\def\FRE{\mbox{Fr\'{e}chet }}

\def\wH{ { \mathfrak{H}}}
\def\TwH{{\mathfrak{H}_\star}}

\def\II{\mathfrak{I}_\beta}

\def\MMM{\mathcal{M}( \vk Y)}	
\def\MMMP{\mathcal{M}_\EI( \vk Y)}

\newcommand{\norm}[1]{\lVert  #1 \rVert }

\definecolor{c20}{rgb}{0.,0.7,0.}
\definecolor{c30}{rgb}{0.,0.,1.}
\definecolor{c40}{rgb}{1,0.1,0.7}
\definecolor{c50}{rgb}{1,0,0}
\definecolor{c60}{rgb}{1,0.9,0.1}
\definecolor{c70}{rgb}{0.50,1.00,0.00}

\newcommand{\Rset}{\mathbb{R}}

\newcommandx{\Norm}[1]{\left|#1\right|}
\newcommandx{\normE}[1]{|\!|#1|\!|}
\newcommandx{\lpnorm}[3][1=,3=]{\left\|#2\right\|_{#1}^{#3}}
\newcommandx{\supnorm}[3][1=,3=]{\left\|#2\right\|_{#1}^{#3}}

\def\boldsymbol#1{#1}
  \def\cadlag{c\`adl\`ag}

\newcommand{\prooftheo}[1]{\underline{Proof of \Cref{#1}:}}
\newcommand{\proofprop}[1]{\underline{Proof of \Cref{#1}:}}
\newcommand{\prooflem}[1]{\underline{Proof of \Cref{#1}:}} 
\newcommand{\proofkorr}[1]{\underline{Proof of \Cref{#1}:}}

\def\R{\Rset}

\makeatletter
\newcommand*\bigcdot{\mathpalette\bigcdot@{.55}}
\newcommand*\bigcdot@[2]{\mathbin{\vcenter{\hbox{\scalebox{#2}{$\m@th#1\bullet$}}}}}
\makeatother

\def \SZ{ \mathcal{S}( {\vk Z})}

\def\SZT{ \mathcal{S}( \ZT)}
\def\ST{ \mathcal{S}( {\vk \Theta})}
\def\SL{ \mathcal{S}( {\vk Q})}

\def\TT{\mathcal{T}}
\def\TTT{\TT}
\def\intT{\int_{\TT}}
\def\TT{\mathcal{T} }

\def\ZT{\widetilde{\vk Z}}   
\def\EIE{\mathcal{L}}
\def\EI{\mathcal{L}}

\def\STYq{\mathcal{S}_\EI(\vk Y)}

\def\kK{\mathscr{C}_{\kappa} [\vk Z]}

\def\kkN{ \mathscr{C}_\kappa[ \vk Z_N]} 

\def\diad{{\mathbb{T}_0}}

\def\wHO{ \wH_{0}}

\def\cadlag{c\`adl\`ag}

\def\bqny#1{ \begin{eqnarray*} #1 \end{eqnarray*}}
\def\bqn#1{ \begin{eqnarray} #1 \end{eqnarray}}

\newcommand{\BEX}{\begin{example}}
\newcommand{\EEX}{\end{example}}

\newcommand{\BIT}{\begin{itemize}}
\newcommand{\EIT}{\end{itemize}}
\newcommand{\BDI}{\begin{description}}
\newcommand{\EDI}{\end{description}}

\newcommand{\nelem}[1]{\Cref{#1}}

\def\mj{\mathcal{J}}

\def\TT{\mathcal{T}}

\def\Z{\mathbb{Z}}
\def\inn{\in \mathbb{N}}

\def\rE#1{\textcolor{c30}{#1}}
\def\rE#1{{#1}}

\def\drz{\lambda_\alpha(dz)}
\def\drzz{\alpha z^{\alpha- 1} \lambda(dz)}
\def\Sset{\mathfrak{V}}

\def\pkZ{p^>_{Z_\kappa}}

\numberwithin{equation}{section}
\def\JF{\mathcal{J}^\star}
\def\JFA{\mathcal{J}^\star_1}
\def\JFB{\mathcal{J}^\star_2}
\def\JFC{\mathcal{J}^\star_3}

\begin{document}

\title{Cluster Random Fields and Random-Shift Representations}

	\author{Enkelejd  Hashorva}
	\address{Enkelejd Hashorva, HEC Lausanne, 
		University of Lausanne,\\
		Chamberonne, 1015 Lausanne, Switzerland}
	\email{Enkelejd.Hashorva@unil.ch}

	\bigskip
	
	\date{\today}
	\maketitle

	\begin{quote} 
		{\bf Abstract:} Cluster random fields (CRFs) play a crucial role in the study of extremes of stationary regularly varying random fields (RFs). In particular, they appear in the {\it  Rosi\'nski representation}  of max-stable and  $\alpha$-stable RFs. In this contribution we introduce CRFs 
		in an abstract setting proving that they are crucial for the construction of {\it shift-generated} classes of $\alpha$-homogeneous RFs. 
		Further,  we investigate the relations between  CRFs, {\it tail RFs} and {\it spectral tail RFs}. Applications discussed in this contribution include new representations of extremal functional indices  and purely dissipative max-stable RFs. 
	\end{quote}

    

\section{Introduction}
 As shown initially in \cite{Davis,BojanS}, 
  CRFs play a fundamental role in the study of extremes of stationary regularly time series. In the literature they appear with different labels. For instance 
in  \cite[Def 5.4.6]{kulik:soulier:2020} and \cite[Def 2.3]{MR4280158} they are referred to as conditional {\it spectral tail RFs}, while 
 \cite{BP, Planic} used the term {\it anchored tail processes}. Adopting the parlance of  
 \cite{BojanPhilippe, Hrovje, kulik:soulier:2020, PH2020,MR4280158, kulik2023,chen2023asymptotic} where the {cluster measures} are defined via CRFs, we shall adhere to the CRF terminology throughout this paper. \\
The first systematic analysis of CRFs in the context of extremes of regularly varying stationary time series (in the discrete-time setting) appeared in \cite{Hrovje}, and was further developed in \cite{klem, kulik:soulier:2020, PH2020, Planic}. 
Crf’s—often also called random shape functions, see e.g., \cite{kab2009,KabExt}—have independently emerged 
in  the {\it Rosi\'nski representation} of max-stable and $\alpha$-stable RFs. In this context, and in connection with {\it Pickands constants}, CRFs constructions for max-stable processes are studied in \cite{Htilt, debicki2017approximation, MR3745388}, as well as in \cite{Hrovje, kulik:soulier:2020, PH2020}.
 
In applications, notably those surveyed in the state-of-the-art monograph  \cite{kulik:soulier:2020}, 
  CRFs are ubiquitous in statistical modelling of stationary regularly varying time series. In particular, CRFs  are the key to the  $m$-approximation technique developed in \cite{klem, kulik:soulier:2020, PH2020}, being also pivotal in the estimation of functional indices, see  \cite{PH2020,MR4280158,kulik2023,chen2023asymptotic}. 
  They are particularly central to the $m$-approximation technique developed in \cite{klem, kulik:soulier:2020, PH2020}, and are instrumental in the estimation of functional indices (see \cite{PH2020, MR4280158, kulik2023, chen2023asymptotic}). Recent developments in \cite{kulik:soulier:2020, MR4280158, chen2023asymptotic} emphasize the use of CRFs in constructing shift representations of {\it tail measures}  via cluster measures (see also \cite{klem, PH2020, kulik2023, wangExt, chen2023asymptotic}).
  
In the context of regularly varying time series, CRFs have been mainly studied in the discrete setting see e.g.,  \cite{BojanS, BojanPhilippe, Hrovje,kulik:soulier:2020,BP,  Planic}. 
An exception is \cite{PH2020}, where cluster processes  with \cadlag\ sample paths and corresponding {cluster measures} are introduced. 

 To set up our mathematical framework, let us fix two positive integers $d$ and $l$, and consider below the parameter  set $\TT=\R^l$ or $\TT= \Z^l$. Fix next $\alpha>0$ and 
write $\Dset$ for  the space of functions $f: \TT\mapsto \R^d$  equipped with  the product (cylindrical) $\sigma$-field   $\AA$. Taking $\mathcal E=\R$ or $\mathcal E=[0,\IF]$,  let $\mathcal{H}$ consist of all  $\AA/\borel{(-\IF,\IF]}$-measurable maps $F: \Dset  \mapsto \mathcal E$, which are bounded if  $\mathcal E=\R$. Let $\mathcal{H}_\beta,  \beta\ge 0$  consists of   all 
$F\in \mathcal{H}$ satisfying   $F(c f) = c^\beta F( f)$ for all $f \in \Dset$ and $c>  0$ and 
set 
$$\mathcal{H}_\beta^+=\{F\in \mathcal{H}_\beta, F\ge 0\}.$$ 
In our notation $\borel{S}$ stands for the Borel $\sigma$-field of a topological space $S$. Throughout the paper $\kappa \in \mathcal{H}_\alpha^+$ is fixed and for a given RF $V(t),t\in \TT$ we write 
\bqn{ \label{vkappa}
	V_\kappa(t)=\kappa(B^{-t} V), \quad B^{-t} f=f(\cF{\cdot+t}), \quad f\in \Dset, \quad t\in \TT.
}
\BD 
Let $\clasP$ be  the class of all $\R^d$-valued RFs $\vk V(t),t\in \TT$ defined on some complete non-atomic probability space such that $V_\kappa$ is stochastically continuous. 
\ED 
Necessary and sufficient conditions for $V_\kappa$ to have a separable and jointly measurable version are given in \cite[Prop 9.4.2]{MR1280932}, which can be formulated with respect to 2-dimensional marginal distributions of $V_\kappa$. In particular,  if $V_\kappa$ is stochastically continuous, then in view of \cite[Thm 5, p.\ 169, Thm 1, p.\ 171]{MR636254}   it has a separable version with separant $\diad$   being further jointly measurable.\\
 Hereafter, elements of $\clasP$ are assumed to be separable with separant $\diad$   and jointly measurable.
 Further $\lambda(\cdot)$  stands for the Lebesgue measure if $\TT=\R^l$ or the counting measure on $\TT$ when the latter is discrete.\\
  Below $\diad$  consists 
of all $t\in \R^l\cap \TT$, which have rational coordinates and set 
$$
\sup_{t\in K} f(t)= \sup_{t\in K \cap \diad} f(t), \quad f\in \Dset, \quad K\subset \TT.
$$ 

\BD Given $\kappa \in \mathcal{H}_\alpha^+$ we call  $\vk Q \in \clasP$ a CRF if 
\bqn{
	\label{dallandysh}
	\pk*{\sup_{t \in \TT} Q_\kappa(t)     >0}=1, 	\quad 	 \quad 	\int_{ \TT}  
	\E*{ \sup_{ t \in [-c,c]^l } Q_\kappa(t-v) } \lambda(dv)< \IF, \quad \forall c>0. 	
}	
\label{DefC}
\ED
The second condition in \eqref{DefC} means that we are interested in locally bounded CRFs.
It turns out that $\sup_{t\in \TT} Q_\kappa(t)$ has a finite expectation, see \eqref{supQ} below.   

\begin{example} \label{examp:kappa}
Hereafter $\norm{\cdot}$ denotes a norm on $\R^d$.  Three natural choices for $\kappa$ are:
\begin{enumerate}[(i)]
	\item $\kappa(f)= \norm{f(0)}^\alpha, f\in \Dset$;
	\item \label{examp:kappa3} $\kappa(f)= (\sum_{i=1}^d \abs{f_i(0)}^\alpha/d)^{1/\alpha}, f=(f_1 \ldot f_d)\in \Dset$;
	\item $\kappa(f)= \sup_{t\in K\cap \diad} \norm{f(t)}^\alpha, f\in \Dset$ for some compact set $K\subset \R^l$.
\end{enumerate}
Given $\kappa\in \mathcal{H}_\alpha^+$  as above,  if 
$L:\R^l \to \R^d$ is a deterministic function such that $\kappa(B^t L), t\inr^l$ is a \cadlag\ pdf and $\intT \sup_{t\in [-c,c]^l} \kappa(B^v L) \lambda(dv)< \IF$ for all $c>0$, then a non-random CRF $Q$ 
is simply 
$
\vk Q(t)=  L(t), t\in \TT.$ 
\end{example}
Consider  $\vk Z \in \clasP$  
defined on a complete non-atomic probability space $(\Omega, \mathscr{F}, \mathbb{P})$ 
 satisfying (recall our notation \eqref{vkappa})
\bqn{\label{dense}
\pk*{\sup_{t\in \TTT} Z_\kappa(t) >0}=1,    \quad 
	\E*{ \sup_{t\in [-c,c]^l } Z_\kappa(t ) }\in (0,\IF), \quad \forall c\in [0,\IF).
} 
Similarly,     $\ZT \in \clasP $ is defined on a complete non-atomic probability space $( \widetilde{\Omega}, \widetilde{\mathscr{F}}, \widetilde{\mathbb{P}})$. 
 We introduce next important classes of $\alpha$-homogeneous RFs, which in the case of {\it shift-generated} dissipative classes (defined below) are directly constructed by CRFs.  
\BD 
\label{xhonga} 
As in \cite{hashorva2021shiftinvariant},  we call $\kK$ with  $\kappa\in \mathcal{H}_\alpha^+$ an $\alpha$-homogeneous class of RFs with representer  $\vk Z$, if it contains   $Z$ and  all  $\ZT \in \clasP$ 
that satisfy  \eqref{dense}  and moreover 
\bqn{ \label{boll} 
\E{ F( \vk Z)}=	\ETT{F( \ZT )}, 
\quad   \forall F\in \mathcal{H}_\alpha.
}
$\kK$ is called {\it shift-generated} (and then denoted by $\kk$) if further 
$$B^h \ZT\in \kK, \quad \forall h\in \TTT$$
for some (and then for all) $\ZT \in \kK$.
\ED 
If the random variable (rv) $C$ is almost surely (a.s.) positive, with $\E{C}=1$ being further  independent of $Z$, then clearly $\ZT= C^{1/\alpha} Z \in \kK$ and 
\bqn{\label{cA}
 \kK= \mathscr{C}_{\kappa} [C^{1/\alpha}\vk Z].
}
Note that in \Cref{xhonga} we require that 
\bqn{\label{PKZ}
	\pkZ=\pk*{\sup_{t\in \TTT} Z_\kappa(t) >0}=1.
}
When $C=0$ with non-zero probability, then $\ZT=C^{1/\alpha}Z$ does not satisfy the first condition of \eqref{dense}. Hence by \Cref{kor21}, \Cref{21:b} below we conclude that $\ZT\not \in \kK$.

In view of \cite{MolchanovBE}, if  $\alpha=1$ and $Z$ is symmetric (non-negative), $\kK$ contains all symmetric (non-negative) RFs $\ZT \in \clasP$ that are zonoid-equivalent to  $Z$. 
Note that therein $( \widetilde{\Omega}, \widetilde{\mathscr{F}}, \widetilde{\mathbb{P}})$ is assumed to be a Borel probability space. The next example  explains that $\kK$'s naturally arise in connection with max-stable and $\alpha$-stable RFs, a fact which is known from \cite{MolchanovBE}.

\begin{example} 
	\label{luxe}
Given   $\ZT \in \kK$ let 
  $ \ZT^{(i)}, i\inn$ be  independent copies of $\ZT$ and define the  max-stable stationary RF $ \widetilde{X} (t),t\in \TT$  via its de Haan representation   (see e.g., \cite{deHaan,dom2016}) by
\bqn{\label{eq1}
	\widetilde{X} (t)=  \max_{i\ge 1} \Bigl(\sum_{k=1}^i \mathcal{V}_k\Bigr)^{-1/\alpha}  \ZT^{(i)}(t), \quad t\in \TT. 
} 
Here the   maximum is applied component-wise, 
with $\mathcal{V}_k, k\ge 1$ unit iid exponential rvs being independent of any other random element.  The RF $\ZT$ is referred to as the representer of $\widetilde X $. \\
For $\alpha \in (0,2)$ assuming further that the elements of $\kK$ are symmetric if $\alpha \in [1,2)$, 
in view of \cite[Lem C1]{PH2020} we can define an $\alpha$-stable RF $\widetilde{X}_{\Sigma}$ by the following  LePage representation
$$\widetilde X_{\Sigma}(t)	=  \sum_{i\ge 1} \Bigl(\sum_{k=1}^i \mathcal{V}_k\Bigr) ^{-1/\alpha}  \ZT^{(i)}(t), \quad t\in \TT. $$
By \cite[Thm 9]{MolchanovBE}, when  $\alpha=1$ the laws of $\widetilde   X $ and $\widetilde X_\Sigma$ do not depend on the choice of $\ZT$. See also \cite[Thm 1.4.2]{MR1280932} covering the case  $\alpha\not=1$. \\
Conversely, if  $\kappa$ is as in \Cref{examp:kappa}, 
given a stochastically continuous max-stable RF $\widetilde X $ with de Haan representation \eqref{eq1} and a representer $Z$ with non-negative components satisfying 
\eqref{dense}, then $\kK$ contains all non-negative RFs $\ZT \in \clasP$ that are valid representers
 for $\widetilde X $ and satisfy $\pk{\sup_{t\in \TT}\kappa(B^t \ZT)>0}=1$. 
\end{example}
Surprisingly, {\it shift-generated} $\kk$'s can be constructed even when $Z$ is not stationary. Their definition 
 is motivated by that of {\it zonoid-stationarity} in \cite{MolchanovBE}  and the characterisation of stationary max-stable RFs in \cite[Thm 6.9]{Htilt}, see also \cite[Eq.\ (5.2)]{Hrovje} and \Cref{sec:maxstable} below. \\
 Hereafter,  the $\TT$-valued rv $N$   with positive  pdf $p_N(t)>0, t\in \TT$ 
is assumed to be independent of any other random element defined in the same probability space. \\ 
It turns out that the construction of  {\it shift-generated} $\kK$'s is closely related to the existence of CRFs.
 Namely, if $Q$ is a given CRF, then letting 
\bqn{ 
	\label{mm} 
	\vk Z_N(t) = 
	\frac{ B^N \vk Q(t)}{ [p_N(N)]^{1/\alpha}  }, \quad t\in \TT,
}
which  in view of  \Cref{molla} is well-defined and belongs to $\clasP$, we obtain that
$\kkN $ denoted simply below by $\kkQ$  is a {\it shift-generated} $\alpha$-homogeneous class of RFs.\\

 By our assumption the RF $Z_\kappa$ is jointly measurable and non-negative. Hence the completeness of  $\Pspace$, the non-negativity of the map $\kappa$ and  the Tonelli Theorem yield 
$$	
\SZ=\intT \kappa(B^{-t}Z)  \lambda(dt)
=\intT Z_\kappa( t)  \lambda(dt)
$$
is a well-defined non-negative rv.
\BD 
If $\pk{\SZ=\IF}=1/0$, then we call $\kk$ purely conservative/dissipative.
\ED 
\BEL 
Let $\kk$ be shift invariant. It is purely conservative/dissipative if and only if (iff) the max-stable RF $X_\kappa$ 
with representer $Z^{1/\alpha}_\kappa$ and de Haan representation as in \eqref{eq1} 
is purely conservative/dissipative. 
\label{kecat} 
\EEL 
Next, we present an important instance of purely conservative $\kK$'s demonstrating in particular their existence. 
\BEL If $\vk Z\in \clasP$ is stationary and satisfies  \eqref{dense}, 
then $\kK$ is {\it shift-generated} and purely conservative.
\label{springWaltz}
\EEL 
If $Z$ is as in \Cref{luxe}, then $\kk$ being purely dissipative is equivalent with $\widetilde X_\kappa $ being also purely dissipative. Conditions for conservativity/dissipativity of $\alpha$-stable RFs are well-known, see e.g.,  \cite{RZ,Genna04,Genna04c,Roy1}. 

Our focus in this paper is on purely dissipative classes of RFs, which in some cases  can also be directly 
constructed by a given $Z$ without any reference to a particular CRF $Q$.\\
 An interesting instance is that of  the  Brown-Resnick $\kk$ introduced in \cite{hashorva2021shiftinvariant}.  The  class of the Brown-Resnick max-stable RFs is discussed in \cite{kab2009,kab2009a,kabluchko2011extremes,kab2010,WangStoev,MR3299407,DombryE,Das,Hrovje, kulik:soulier:2020, KumeE} motivated by the two prominent instances of stationary max-stable processes presented  in \cite{bro1977,eddy1980distribution,Gale}. 

\begin{example} (Brown-Resnick $\kk$)
Let $\kappa$ be as in \Cref{examp:kappa}	 and define
	\bqn{ \label{BR}
		\vk Z(t)= (   e^{    W_1(t) -  \E{ W_1^2(t)}/2} \ldot 
		e^{    W_d(t) -  \E{  W_d^2(t) }/2} ),     \quad  t\in \TT, 
	}
	where $\vk W(t), t\in \TT$ is a centered $\R^d$-valued Gaussian RF with stochastically continuous sample paths.   Define the matrix-valued pseudo-cross variogram 
	$\gamma$ by  
	$$\gamma_{ij}(s,t)=\E{ (W_i(t)- W_j(s))^2}, \quad 1 \le i,j\le d, s,t\in \TT$$
	If we suppose further that $\gamma_{ij}(s,t)$'s depend  only on $t-s$ for all $s,t \in \TT$ and all 
	positive integers $i,j\le d$, as shown in \cite[Example 4.2]{hashorva2021shiftinvariant}, 
	see also \cite[Lem 4.2]{KumeE} and  \cite{MolchanovBE,oesting2017statistical},  it follows that $\kK$ is shift-generated. 
In view of  \cite{DombKabB},  see also \cite[Rem 15]{kab2009}, setting below 
	$$
	\norm{h}_*=\sum_{i=1}^l \abs{t_i}, \quad  h=(h_1 \ldot h_l)\in \R^l
	$$
and  using \Cref{kecat}, we conclude that $\kk$ is purely dissipative, provided that 
\bqn{\label{BR2}
		\liminf_{ h\to \IF }  \max_{1 \le i\le d}\frac{ \gamma_{ii}(0,h) }{ \ln \norm{h}_*}> c
}
for some $c>0$ sufficiently large. \\
\rE{If  $\kK$ is  shift-generated, 
	the law of the corresponding max-stable RF $X$   depends only on the cross variogram. Moreover, if  a corresponding CRF $Q$ can be constructed, then also the law of $Q$ depends only on $\gamma$.}
\label{ex:BRModel}
\end{example}

Note in passing that the previous example indicates that not every   $\vk Z\in \clasP$ that satisfies \eqref{dense} defines a {\it shift-generated} $\alpha$-homogeneous class of RFs.  One such instance is $\vk Z$ in \eqref{BR} for which $\gamma_{ij}(s,t)$'s depends on both $s,t$ and not only on the difference $t-s$. \\
We consider next the Brown-L\'evy-Resnick $\kk$ discussed in \cite{hashorva2021shiftinvariant}. In the study of max-stable  processes it has been introduced in \cite{StoevSPA}, see \cite{eng2014d,Htilt} for further results.
\begin{example}(Brown-L\'evy-Resnick $\kk$) 
	\label{exa:lev} Let $W_i(t_i),t_i\ge 0,i\le l$ be L\'evy processes with Laplace exponent $\psi_i(\theta)=\ln \E{ e^{\theta W_i(1)}}$ such that $\psi_i(\alpha)=0$ for some $\alpha>0, i\le l$. Write $W_i^{(\alpha)}(t_i),t_i\ge 0$ for the L\'evy process with Laplace exponent $\psi_i(\alpha+\theta)$. Assume that $W_i, W_i^{(\alpha)}$'s are all independent with \cadlag\ sample paths and define 
	$$ Z(t)= \prod_{i=1}^l e^{ \ind{t_i\ge 0}W_i(t_i)  - \ind{t_i<0}  W_i^{(\alpha)}((-t_i)-)}, \quad t=(t_1 \ldot t_l).
	$$ 
	It follows that the max-stable RF $X $ defined from $Z$ 
	via \eqref{eq1} is stationary (for $l=1$ see e.g., \cite{eng2014d}). Since further $Z$ 
	satisfies \eqref{dense}, then $\kK$ is shift-generated. Using that both $W_i(t_i), W_i^{(\alpha)}(t_i)$ converge to $-\IF $ as $\abs{t_i}\to \IF$, then in view of \Cref{kecat} and \Cref{kabA} we conclude  that $\kk$ is purely dissipative.
\end{example}

\BD 
We say that $\kK$ has a   random-shift representation, if  
$\vk Z_N$ defined in \eqref{mm} for some $Q \in \clasP$ belongs to $\kK$. 
\ED  
When $\kK$ has a {\it random-shift representation},  then it is necessarily {\it shift-generated} and 
{\it purely dissipative}, see \Cref{conv} below. Moreover, it agrees with $\kkQ $ and necessarily   $\vk Q$ is a  CRF.\\

In both examples considered above it is natural to ask the following question:
\\
Does  the {\it shift-generated} $\kk$ has a  {\it random-shift representation}  and if so, how to determine a corresponding CRF $Q$?\\
 An equivalent question concerns the existence of the {\it  Rosi\'nski representation}  (according to the terminology of \cite{Roy1}). Namely, it is of interest if a \cadlag\ RF $Q(t),t\in \TT$ exists 
  such that 
\bqn{\label{eqcR} 
	\widetilde X (t)= \max_{i\ge 1} P_i Q^{(i)}(t-T_i), \quad t\in \TT,
}	
 where $\{P_i, T_i \}'s$ are points of a Poisson Point Process on $(0,\IF) \times \R^l$ with mean measure 
$ \lambda_\alpha(\cdot)  \odot \lambda(\cdot)  $ and $Q^{(i)}$'s are 
independent copies of $Q$ being further independent of 
$\{(P_i, T_i), i\ge 1\}$; throughout the paper $ \lambda_\alpha(dr)= \alpha r^{-\alpha-1} d r$.\\
	For  max-stable Brown-Resnick RFs the existence of a corresponding CRF $Q$ has been studied in  \cite{kab2009, MR3215542,dom2016}. The law of $Q$ is however in general not known.  For specific cases, e.g., $\TT=\Z^l$ it has been determined in \cite[Thm 8]{KabExt}. \\
Utilising another approach, CRFs for \cadlag\ max-stable RFs 
 are constructed in \cite{debicki2017approximation}. See also \cite{klem,kulik:soulier:2020,PH2020} for related results and ideas.\\
The connection of CRFs with {\it  Rosi\'nski representation}  of max-stable RFs is first shown in \cite[Thm 5.1]{Hrovje} (the construction is similar to that of  \cite[Thm 8]{KabExt}) and then further discussed in \cite{kulik:soulier:2020,klem,PH2020}.
Such representations have initially appeared in fact earlier,  for instance 
in the study of $\alpha$-stable processes,  see e.g., \cite{RZ,RZ2,Genna04,Genna04c,kab2009a, WangStoev,MolchanovBE, Roy,Roy1,PH2020}.\\
Summarising the findings in literature, 
the principal applications of CRFs concern: 
\begin{enumerate}[i)]
	\item Explicit construction of 
	{\it cluster measures}  discussed in \cite{kulik:soulier:2020,PH2020} as well as determination of 
	{\it tail and spectral tail RFs} introduced in \cite{BojanS,BP}. 
	\item Estimation of functional indices, see e.g., \cite{kulik:soulier:2020,PH2020,MR4280158, kulik2023}.
	\item Representations  of  {\it Pickands constants} and 
	extremal functional indices, see \cite{Htilt,HBernulli,debicki2017approximation,PH2020,BP,Planic,hashorva2021shiftinvariant}.
	\item $m$-approximation of stationary regularly varying RFs, see e.g., \cite{kulik:soulier:2020,BojanPhilippe,klem}.
	\end{enumerate} 
Natural questions that arise in the general setup of this paper, i.e., dropping the \cadlag\ assumption on the sample paths,  include:   
\begin{enumerate}[Q1)]
	\item \label{refQA}
Given a {\it shift-generated} $\kk$, under what conditions on $\vk Z$ does it possess a  {\it random-shift representation}  and how to determine a corresponding CRF  $\vk Q$? 
	\item \label{refQC} How  to construct different CRFs $\vk Q$?
	\item \label{refQE} How do different CRFs $Q$ define the same {\it Pickands constants}?
	\item \label{AppBR} \rE{Is the class of {\it shift-generated} Brown-Resnick $\kk$'s only dependent on $\kappa$ and $\gamma$ and what about the law of $Q$?}
\end{enumerate}
Answering \Cref{refQA} is important since then Rosi\'nski representations of max-stable and $\alpha$-stable RFs can be easily addressed, see \Cref{sec:maxstable} where we discuss in particular new representations for both the Brown-Resnick and the Brown-L\'evy-Resnick max-stable RFs.\\
 A key finding used in our constructions is \cite[Thm 2.1]{GeoSS}. 
 Additionally, since the integral functional is not measurable for the product $\sigma$-field, another crucial  result needed in the proofs is the extension of \eqref{boll}, presented below in \Cref{thmAug23}. 

\Cref{refQC} leads to different Rosi\'nski representations, which in turn imply new tractable expressions for extremal functional indices, see \Cref{sec:candi}. In applications, having different representations for the extremal index is important since it allows for construction of flexible estimators. If one is interested in calculating those indices, as in the case of {\it Pickands constants} , such representations can be further  useful for both Monte Carlo simulations and derivations of precise bounds.  We note in passing that the classical Pickands constant is not known apart from two particular values and its calculation is still an interesting research topic, see e.g., \cite{DiekerY, bisewski2021speed,adonai} and the references therein.\\
 In our study of CRFs we obtained also new results for {\it shift-generated} $\kk$'s as well. For instance, the identity presented in \eqref{inemato} extends a previous one derived in \cite{DiekerY} for $k=d=1$.

Brief outline of the rest of the paper:  We shall present  the main notation and definitions in \Cref{sec:prim}, which includes also few preliminary results. In particular, \Cref{kabA} derives new equivalent conditions, which are important for the determination of purely dissipative $\kk$'s.\\ 
\Cref{MainRes} answers \Cref{refQA} and \Cref{refQC} by discussing first   
 basic properties of CRFs followed by explicit constructions of  {\it random-shift representations}  for $\kk$ based on results and ideas presented in \cite{klem,Hrovje,kulik:soulier:2020,PH2020,hashorva2021shiftinvariant}. 
All the  constructions in \Cref{MainRes}  are new  if $\norm{\vk Z(0)}$ is a.s.\ positive, which is in particular the case for both the Brown-Resnick and the Brown-L\'evy-Resnick $\kk$'s. 
\Cref{sec:cash} discusses several applications and answers in particular \Cref{refQE} and \Cref{AppBR}.  Concluding, in \Cref{appn} we present some technical results.

\section{Preliminaries}
\label{sec:prim}
We shall introduce first some classes of maps. In  \Cref{sotoS} {\it local RFs}, {\it spectral tail} and {\it tail RFs} are briefly discussed followed by a short investigation of conditions that characterise pure conservativity/dissipativity in \Cref{adonaiPD}.

\subsection{Homogeneous, anchoring and involution maps}
\label{some:class}
Hereafter  $H : \Dset \to [-\IF, \IF]$  is  called shift-invariant if 
(recall the definition of $\mathcal{H}$ and $\mathcal{H}_\beta$ in the Introduction) 
\bqn{ \label{mapShift}
	H(B^hf) =H(f), \quad \forall f\in \Dset,h\in \TT.
	}
 
When $\TT=\R^l$, for common choices of $\kappa$  the following integral map
$$
F_I(f):f \mapsto \intT \kappa( B^{-t} f)^\xi\lambda(dt),\quad f\in \Dset, \quad \xi >0 
$$	
is  not $\AA/ \borel{[0,\IF]}$-measurable (set  $F_I(f)=+\IF$ if it is not defined).  \\ 
\cF{Indeed, if for instance $\kappa$ is specified as in \Cref{examp:kappa}, the non-empty set $A=\{f\in \Dset: F_I(f)= 0\}$ consists of $f$'s that vanish almost everywhere on $\TT$ and clearly $A \not=\Dset$. If $A\in \AA$, since the elements of the product $\sigma$-field $\AA$ on $\Dset$ have only countable restrictions (see e.g., [Lem 1.5. 52]) that  are irrelevant for $F_I(f)$ (recall that $\lambda(\cdot)$ is the Lebesgue measure on $T$), then $A=\Dset$, which is a contradiction.  
}	

\def\II{\mathfrak{I}}
Let $g_i:\TTT\mapsto [0,\IF), i\le 3$ be locally bounded and  $\lambda$-measurable. When $\TT=\R^l$ it shall be assumed that $g_i$'s are positive  almost everywhere. Write  $\ind{A}$ for   the indicator function of some set $A$. 

\BD  \label{twq}
 Write   $\Hh_\beta, \beta \ge 0$ for the subset of $\mathcal{H}_\beta$ \cED{with elements $F$} 
determined by 
\bqn{
\label{eier}
	F(f)= \frac{ \Gamma (f) \II_1  (f) \ind{\II_2(f) \in A}  }{ \II_3  (f)  }, \ \  \ \II_i(f) = \int_{\TTT} \kappa(B^{-t}f)^{\xi_i}  g_i(t) \lambda(dt), \quad f \in \Dset,
}
where $\Gamma  \in \mathcal{H}_{\xi_0}, g_i, i\le 3$ with  constants  
	$\xi_i$'s  and  
$ A \subset [0,\IF]$ a Borel set such that $F$ is $\beta$-homogeneous. 
 If $F(f)$ is undefined we set   $F(f)=+\IF$. 
Write $\TwH$ for the class of maps $\Gamma $ defined for some $\beta\in [0,\IF) $ by 
$$
\Gamma= F_1F_2, \quad   F_1 \in\Hh_{\beta}, F_2\in \mathcal{H} .
$$ 
\ED
It is possible to extend \eqref{boll} to include $F$ in $\Hh_\beta$. The corresponding findings of \cite{hashorva2021shiftinvariant} are summarised next: \BEL For all $Z \in \clasP, \Gamma\in \TwH$ the law of $\Gamma(\ZT)$ 
depends only on the finite dimensional distributions (fidi's) of $\ZT \in \kK$ and further 
\bqn{ \label{tcfN} 
	\E{ F( Z )} = \ETT{F(   \widetilde{\vk Z})}, \quad \forall F\in \Hh_{\alpha} , \forall   \ZT \in \kK.
} 
Moreover, $\kK$ is {\it shift-generated}   iff \eqref{tcfN} holds for all $\ZT=B^h Z, h\in \TT$ and it is also equivalent with
\bqn{ \label{tcfN:0} 
	\E{ \kappa(B^{-h} Z) F(Z)} = 	\ETT{ \kappa(\ZT) F(B^h \ZT)}, 
	\quad \forall	F\in \Hh_{0}, \forall h\in \TT,\forall \ZT \in \kK.
} 
\label{thmAug23}
\EEL  
\BRM 
In view of \eqref{tcfN} 
the shift-generation of $\kk$ is a property of $Z$ and not of $\kappa$. However, by the definition the pure conservativity/dissipativity relates to $\kappa$.
\ERM

\subsection{Local, spectral tail and tail RFs}
\label{sotoS}
Tail and spectral tail RFs play a crucial role in the asymptotic analysis of multivariate regularly varying time series. Initially introduced in   \cite{BojanS}, these RFs have been studied in numerous contributions, see e.g., \cite{SegersEx17, Hrovje,kulik:soulier:2020,BP,MR4280158,MR4280158,MartinE, MolchRV}.\\
In the sequel  $\kappa \in \mathcal{H}_\alpha^+$ is fixed,  
$\kk$ is {\it shift-generated} and  we shall suppose for simplicity that  $\E{\kappa(  Z)}=1$ (recall that we assume \eqref{dense}).  Consequently, in view of  \eqref{tcfN} 
\bqn{\label{perfect}
	\ETT{\kappa(B^{t} \ZT)}=1 , \quad\forall  t\in \TT,  \forall \ZT \in  \kk.
}
Write $\widetilde\Theta$  for the RF  $ \widetilde Z/ \kappa^{1/\alpha}( \widetilde Z) $ under the probability measure 
\bqn{\label{minist}
	\widehat {\mathbb{P}}(A)=  
	\ETT{   \kappa(   \ZT )  \ind{ A}}, 
	\  \
	\quad \forall A \in \mathscr{F}.
} 
Note that since $(  \widetilde{\Omega},  \widetilde{\mathcal{F}},  \widetilde{\mathbb{P}})$  is complete and non-atomic, then also 
$(  \widetilde{\Omega},   \widetilde{\mathcal{F}},\widehat{\mathbb{P}})$ is complete and non-atomic.  Let  hereafter  $R$  be  an $\alpha$-Pareto rv with survival function $s^{-\alpha}, s\ge 1$) defined on $ (\widetilde{\Omega}, \widetilde{\mathscr{F}}, \widehat{\mathbb{P}})$   being independent of any other random element. 
\BD 
We shall call $ {\Theta}$ corresponding to $Z$ a local RF  of $\kk$ 
and set  $ {Y}(t)=R  {\Theta}(t),t\in \TT$.
\ED 
\rE{\BRM\label{thetachanges} 
The local RF $\Theta$ of a given $\kk$ changes in general,
 if we choose to tilt with respect to some $\kappa_\star \in \mathcal{H}_\alpha^+$. In particular, if 
		\bqn{ \kappa_\star(f)>0 \implies \kappa(f)>0, \quad \forall f\in \Dset,
				\label{simhusimhu}
			}
then $\Theta_\star$ that corresponds to the tilted law of $Z/ \kappa_\star^{1/\alpha}(Z)$ satisfies 
\bqny{
	\EHh{F(\Theta_\star)}= \ETT{ \kappa_\star( Z) F(Z/\kappa_\star^{1/\alpha} (Z)) }
= 
	\EHh{  \kappa_\star( \Theta )  F(\Theta/ \kappa_\star^{1/\alpha}(\Theta ) )}, \quad \forall F \in \mathcal{H}.		
}
Consequently, the  law of $\Theta_\star$  is determined by $\Theta$.\\
\ERM
}

As shown in \cite{hashorva2021shiftinvariant}, see  \cite{PH2020} for the \cadlag\ case and $\tau=0$, the second condition in \eqref{dense} is equivalent with 
\bqn{ \label{leviz}   
	\EHh*{ \frac{ 1}{\int_{  [-c,c]^l \cap \TT }  \kappa^\tau( B^{t-s}  \Theta ) \ind{ \kappa(B^{t-s} \vk Y) > 1 }  \lambda(ds)}} \in (0,  \IF)
}
for all $c>0$ and all $\tau\inr$ such that 
\bqn{ \label{tauR} 
	\sup_{s\in [-c,c]^l \cap \TT}\EHh{\kappa^\tau( B^{-s} \Theta)  }< \IF, \quad \forall c>0.
}

In view of \cite{hashorva2021shiftinvariant} Eq.\ \eqref{tcfN:0} implies for all $ h\in \TT,x>0$  (recall the definitions of maps in \Cref{twq})
\bqn{\label{eqDo20}
  \EHh*{ \Gamma( B^{h}  {\vk \Theta})\mathbb{I}(  \kappa(B^h    {\vk\Theta} )\not =0) } 
  &=& 	\EHh*{\Gamma({  {\vk \Theta}})\kappa( B^{-h}   {\Theta }) } , \quad 
\forall \Gamma\in  \wHO,\\
x^{-\alpha}\EHh{  \Gamma( xB^h  {\vk Y})\ind{ \kappa( x B^h  {\vk Y} ) >1 }} &=& 
\EHh{  \Gamma(   {\vk Y})  \ind{ \kappa( B^{-h}   {\vk Y } /x ) >1 }   }, \quad   \forall \Gamma  \in \TwH.
\label{tYY}
}  
In \cite{Htilt} the local RFs were shown to satisfy the time-change formula \eqref{eqDo20}.  The functional identity \eqref{tYY} for {\it tail RFs} $Y$ is first derived  in \cite{Hrovje}. As in \cite{hashorva2021shiftinvariant}, we can define a {\it spectral tail rf} $\Theta$ and the corresponding $Y$ directly, without any reference to $Z$, see \cite{PH2020,hashorva2021shiftinvariant} for characterisation results. 

\BRM \label{remMaps} 
\begin{enumerate}[(i)]
\item \label{remTHx:1} In both \eqref{eqDo20} and \eqref{tYY} we interpret $\IF \cdot 0$ and $0/0$ as 0 and those rules apply also hereafter.
\item \label{remMapsP} If $\kk=\mathcal{C}_{\kappa} [\overline{ Z} ] $, 
then  the law of $\Theta$ is the same as the law of the local RF $\overline{ \Theta}$ of $\mathcal{C}_{\kappa} [\overline{ Z} ]$. The converse claim is also valid. This fact explains in particular the advantage of working with the local RF instead of $Z$, since it is unique in law. We note further that for all 
$\Gamma\in \TwH$ the law of $\Gamma(\Theta)$ depends only on fidi's of $\Theta$. 
\item For discrete $\TT$ the equivalence of \eqref{tcfN:0} and \eqref{tYY} is first derived   in \cite{Hrovje}, see \cite{klem, BP,Planic,PH2020,MartinE} 
for  extensions and \cite{BojanS} for the initial formulation of \eqref{eqDo20}.
 	\item  \label{remMaps:2} In the context of regularly varying time series, {\it spectral tail RFs} are defined without reference to $Z$, see \cite{BojanS,kulik:soulier:2020,Resnickart}. Their definition directly from  $Z$ is first shown in \cite{Htilt}, see also \cite{hashorva2021shiftinvariant, LastG}.
\end{enumerate}
\ERM   

\subsection{Characterisation of pure conservativity/dissipativity}
\label{adonaiPD}
Our definition of pure conservativity and dissipativity of a {\it shift-generated} $\kk$ agrees with those of the stationary max-stable RF $X_\kappa(t),t\in \TT$ with 
  representer $Z_\kappa(t)= \kappa(B^{-t} Z),t\in \TT$. \\
  A well-known condition for pure dissipativity is $\pk{\mathcal{S}(Z_\kappa)< \IF}=1$, see e.g., 
\cite{Genna04,PH2020}. Several equivalent condition are obtained in the recent contributions, see e.g., \cite{dom2016,klem,kulik:soulier:2020,PH2020}. 

 \def\EIEA{\mathcal{K}}
 Define next for non-empty $\EIEA\subset \TT$ with  $\lambda( \EIEA)>0$ if $\TT=\R^l$   
\bqn{\label{nowWBack}
	 \mathcal{S}_{\EIEA}(f) = \int_{\EIEA} \kappa(B^{-t} f )  \lambda(dt), \quad   
 \mathfrak{B}_{\EIEA, \tau}(f)=
 \int_{\EIEA}  \kappa^\tau(B^{-t}f) \ind{ \kappa(B^{-t}f)>1}\lambda(dt), \quad  \tau\inr.
}
 
An important property of $\mathcal{S}_\EIEA,  \mathfrak{B}_{\EIEA, \tau} $ is their shift-invariance when $\EIEA$ is an additive  subgroup of $\TT$. 
If $\EIEA=\TTT$ we write simply $\mathcal{S}(f)$ instead of $\mathcal{S}_\EIEA(f)$.

Next, let $\EI $ be a countable subset of $\TT$ with infinite number of elements 
and  write $(\R^k)^*=(\R^k \cup \{\IF\}), k\inn$ for  the one-point compactification of $\R^k$. For several instances, we shall consider maps $F$ being shift-invariant with respect to $\EI$. \\
We introduce below three other maps, which are also shift-invariant (with respect to $\EI$ and not $\diad$),
if further $\EI$ is an additive subgroup of $\TT$. 

\BD  Let    $\mathcal{J}: \Dset   \to  (\R^d)^* $  be $\AA/\borel{(\R^d)^*}$-measurable: 
\begin{enumerate}[{J}1)] 
	\item 	\label{A1}	
\rE{	For all $f\in \Dset $   
	\bqn{
		\label{eqA1}
		\JF[f]:=\sum_{j\in \EI} \ind{ \mathcal{J}(B^j f)=0} \le 1.
	}
}
	\item \label{A2}  For all $ f\in \Dset$ if $\mathcal{J}(f)=j\in \EI$, then  $\kappa(B^{-j}f) \ge \min(\kappa(f),1 )$.

	\item
	\label{A3}  For all $f\in \Dset$ if $\mathcal{J}(f)=j\in \EI$,
	then  $\kappa(B^{-j}f) >0$.
\end{enumerate}	
Suppose that  $\mathcal{J}$  satisfies   \ref{A1}. When  \ref{A2} holds it is referred to as  anchoring. If $\mathcal{J}$ is  0-homogeneous it is called 
a shift-involution and if further \ref{A3} is satisfied,
 it is referred to as  a positive shift-involution. 
\ED 
\cED{Note in passing that \Cref{A1}  follows if 
	\bqn{ \label{lere} 
		\mathcal{J}(B^j f) = \mathcal{J}( f)+j, \quad \forall j \in \EI,\forall f\in \Dset,
	}
	which has been assumed in \cite{BP,kulik:soulier:2020,HBernulli, PH2020} in the definition of anchoring maps when  $\EI=\Z$.}

Hereafter  $\prec $ stands for a given total order on $\TT$, which is shift-invariant, i.e., $i\prec j$ implies $i+k \prec j+k$ for all $i,j,k\in \TT$.  We write $i \preceq j$ if $i\prec j$ or $i=j$.  Below, both  $\inf $ and $\sup $ are taken with respect to  $\prec$ order and the infimum of an empty set is equal to $\IF$.  As in \cite{BP} define the first exceedance functional  $\Ife $  by
$$ \Ife(f)= \inf\Bigl(j\in \EI :  \kappa(B^{j}  f)>1  \Bigr), \quad f \in \Dset,
$$
where  $\Ife(f)=  \IF$ if there are infinitely many exceedance on $\{ j \in \EI, j \prec  k_0 \}$ for some 
$k_0\in \EI$  with  all components positive. 
Define further 
$$\Ifasup(f)= \inf \Bigl( j \in \EI:  \kappa(B^{j}  f)= \sup_{i \in \EI} \kappa(B^{ i}  f) \Bigr), \ \  f\in \Dset ,$$
which is   a positive shift-involution being also anchoring.   
If the infimum is not attained at some element of $\EIE$, then the maps defined above are assigned  $ \IF$.

Hereafter, an important instance  for $\EI$ shall be a discrete subgroup of the additive group $\TT$ with  infinite, but countable number of elements, also referred to as a lattice on $\TT$.   For such a lattice, we can find an  $l\times l$ real matrix $A$ (called a base matrix) such that $\EIE=\{ A x, x\in \Z^l\},$ where  $x$ denotes an $l\times 1$ vector.  Two base matrices $A,B$ generate the same lattice iff $A= BU$, where $U$ is an $l\times l$ real matrix with determinant $\pm 1$. 
Denote 
the  fundamental parallelogram of $\EI$ by 
$$P(\EIE)=\{ A x, x\in [0,1)^l\}.$$ 
The volume of the fundamental parallelogram does not depend on the choice of $A$ and is given by 
$$\Delta(\EIE)= \abs{det(A)}.$$ 
\BD We call $\EIE$ a full rank lattice if $A$ is non-singular. 
\ED
Below all set inclusions or set equalities  are  modulo null sets with respect to $\mathbb{P}, \widehat{\mathbb{P}}$ or $\widetilde{\mathbb{P}}$, depending on the context.  
Next suppose that given $\mathcal{J}_i,i=1,2,3$  (recall the definition in \eqref{eqA1})
\bqn{        
	\{\cED{\JFA[\vk \Theta]=1}\} \subset 
	\{ \mathcal{S}_{\EI}(\vk \Theta) < \IF \} , 
	\quad \{ 
	\cED{\JFB[\vk Y ]=1} \} \subset \{ 
	 \mathfrak{B}_{\EI,\tau}(\vk Y)<\IF\}
\label{peres2}
}
for all $\tau \inr$ satisfying \eqref{tauR} and further
\bqn{  
		\{ \cED{\JFC[\vk Z ]=1} \} \subset \{ \mathcal{S}_\EI(\vk Z) <\IF\}.
\label{peres1}
}	 
A particular instance when the above conditions are satisfied is 
$\mathcal{J}_1=\mathcal{J}_3=\Ifasup, \mathcal{J}_2=\Ife$ (recall  the definition in \eqref{eqA1}).

\def\MMM{\mathcal{M}( \vk Y)}	
\def\MMMP{\mathcal{M}_\EI( \vk Y)}	
\def\BBYt{\mathfrak{B}_{\TTT,\tau}(\vk Y)}  
\def\EI{\EIE} 

\BT  For all $\tau\inr$ satisfying \eqref{tauR} we have 
		\bqn{ 
	\quad 	\pk{\SZ>0}=	\pkh{\ST >0}=\pkh{\BBYt >0}=1 . 
	\label{am0}
}	
Let    $\EI$ be a full rank lattice on $\TT$. If $\mj_1$, $\mj_2$ satisfy \eqref{peres2} and $\mj_3$  satisfies \eqref{peres1}, then for all $b \in [1,\IF)$  
\bqn{ \label{am10}
\{ \SY< \IF \}  & =&    
\Bigl \{ \lim_{ \norm{t}_* \to \IF, t\in \EI } \kappa(B^{-t}\vk \Theta ) =0   \Bigr \}=
\Bigl \{ \lim_{ \norm{t}_* \to \IF, t\in \TT } \kappa(B^{-t}\vk \Theta ) =0   \Bigr \}\\
&=& 
  \Bigl\{ \int_{\TT} \sup_{t \in  [-c,c]^l  \cap \diad  } \kappa( B^{s-t} \Theta) \lambda(ds) < \IF \Bigr\}
\\ 
& =&  
\{ \JFA[\Theta ]=1 \}=\{ \JFB[Y ]= 1\}\label{arm10}\\ 
& =&	
\{\STYq <\IF \}		 = \{\mathfrak{B}_{\EI, \tau}(b \vk Y) <\IF  \}
	= \{\mathfrak{B}_{\TT, \tau}(b \vk Y) <\IF  \}	\label{am12}
}	
and 
\bqn{ \label{am20}
\{  \SZT<\IF \} & =&   	 
\Bigl \{ \lim_{ \norm{t}_* \to \IF, t\in \EI } \kappa(B^{-t} \ZT ) =0   \Bigr \}
	= \Bigl \{ \lim_{ \norm{t}_* \to \IF, t\in \TT } \kappa(B^{-t} \ZT ) =0   \Bigr \}  \\
&=& \Bigl\{ \int_{  \TT} \sup_{t \in   [-c,c]^l  \cap \diad  } \kappa( B^{s-t}\ZT ) \lambda(ds) < \IF \Bigr\}\\
&=&\{ \JFC[\ZT ]= 1\},  
\quad \forall \ZT \in \kk.
\label{am20c}
}	
\label{kabA}
\ET

\BRM  \begin{enumerate}[(i)]
	\item   In case of RFs with \cadlag\ sample paths, the claims in  \eqref{am0}  are direct consequences of  the properties of $\Theta$ and $Y$ (recall $\pk{\kappa(\Theta)=1}=1$). An important result which implies \eqref{am0} in the settings of this paper is obtained in  
	\cite[Thm 2.1]{GeoSS}.
	\item If $\mathcal{J}$ is the infargmax or the first/last  exceedance map, then 
	$ \mathcal{J}^*[f]=1$ is equivalent with $\mathcal{J}(f) \in \EI$ (recall that 
	$\EI$ has infinite but countable number of elements).
	\item For $Z_\kappa$ as in \Cref{springWaltz} we have that $\pk{  \SZ=\IF }=1$ and therefore all events in \Cref{kabA} have probability 0. This is in particular the case if $Z$ and $\kappa$ are as in \Cref{ex:BRModel} and $W$ is stationary with positive variance function. 
\end{enumerate}
\ERM
\BEX \label{BRownE} Consider the {\it shift-generated}  Brown-Resnick $\kk$ introduced in \Cref{ex:BRModel}. In view of \cite[Example 4.2]{hashorva2021shiftinvariant} the law of $\Theta$ depends only on the cross variogram $\gamma$. Suppose therefore without loss of generality that a.s.
$$ W_i(0)=0,\quad i=1 \ldot d, \quad \kappa(Z(0))=1.$$
Hence by the shift-invariance 
$\E{ \kappa(B^tZ )}=1$ for all $t\in \TT$. Consequently, $\Theta$ has the same law as $Z$ 
and moreover 
$$
 Y_i(t)= e^{  \mathcal{E}/\alpha}\Theta(t)= e^{  \mathcal{E}/\alpha+ W_i(t)- \alpha \E{W_i^2(t)}/2}, \quad i=1 \ldot d, \quad t\in \TT,
$$
with $\mathcal{E}$ a unit exponential rv independent of the other random elements  
(note that $e^{ \mathcal{E}/\alpha}$ is an $\alpha$-Pareto rv). \\
If \eqref{BR2} holds, then $\pk{\SZ < \IF}=1$. Hence all the events defined in \Cref{kabA} hold with probability one. In particular a.s. 
\bqn{ 
	\int_{\R} e^{ \alpha W_i(t)- \E{(\alpha W_i(t))^2}/2}dt \in (0,\IF), \quad i=1 \ldot d,
\label{zbalie}
}
which for $\alpha=d=1$ has initially appeared in connection with the dissipativity of the corresponding max-stable process in \cite{WangStoev, dom2016}.
\EEX
If $\kk$ is purely dissipative, then 
\cite[Lem 9.11]{hashorva2021shiftinvariant} and \eqref{dense} imply
\bqn{ \label{adonai} 
\pk*{\sup_{t\in \diad} Z_\kappa(t) \in (0,  \IF)}=1.
	}
Under an additional assumption the next lemma states the converse result. 
\BEL 
For a given $\kk$ if \eqref{adonai} is satisfied and further 
\bqn{
	\label{intKab} 
	\int_{0}^\IF t \pk*{ \sup_{s\in \diad} Z_\kappa(s)  \ge t } dt< \IF,
 }
then $\kk$ is purely dissipative. 
\label{kL}
\EEL
 
\def\intL{\int_{\EI}}
\def\BBYt{\mathfrak{B}_{\TTT,\tau}(\vk Y)}  
\def\MM{\mathcal{M}(\vk Y)}
\def\MME{\mathcal{M}_\EIE(\vk Y)}
\def\MMA{\mathcal{M}_{\EIE \cap \Sset}(\vk Y)}
\def\cEp#1{#1}
\def\EIS{\EIE \cap \Sset}
\def\tQ{\widetilde{Q}}

\section{Main Results}
\label{MainRes}
In the first part of this section  we shall discuss  basic properties of CRFs and their relations with purely dissipative {\it shift-generated} $\alpha$-homogeneous classes of RFs. The second part is dedicated to explicit constructions of CRFs $\vk Q$ in connection with {\it random-shift representations}  of a given $\kk$. 

For simplicity we shall assume \eqref{perfect} in the following. 
\subsection{Shift-generated $\alpha$-homogeneous classes and  CRFs}
Recall that for a given  CRF $Q\in \clasP$ and  $\vk Z_N$ determined in \eqref{mm} we write simply $\kkQ$ instead of $\kkN $. 
\BEL 
Given a $\kk$, for all shift-invariant $F\in \Hh_{0}$ 
\begin{enumerate} [(i)]
	\item \label{21:a}$\EHh{F(\vk \Theta)}=0$ for some  (and then for all) local RFs $   \Theta$;
	\item \label{21:b} $\ETT{F(\ZT)}=0$ for some (and then for all) $\ZT \in \kk$; 
	\item \label{21:c} $\E{ F(\vk Q)}=0$ for some (and then for all) CRF $Q\in \clasP$ such that $\kk=\kkQ$ 
		\end{enumerate}
are all  equivalent, where \Cref{21:c} is valid under the additional assumption that 
 $\kk= \kkQ$. 
Moreover, we have  
\bqn{ \label{e21}
	\pkh{ \ST \in (0,\IF)} =  \pk{ \SZ \in (0,\IF)} =  
	\pk{ \SL \in (0,\IF)}=1.
}
Conversely, if $\pkh{ \ST \in (0,\IF)}=1$, 
then $\kk=\kkQ$ with $\vk Q= c^{1/\alpha}\vk \Theta, c=1/\ST$.
\label{kor21}
\EEL

\BD \label{dfDc} 
$\DDset$ denotes the set of all \cadlag\ functions 
$f: \TT \mapsto \R^d$  if $\TT=\R^l$ and $\DDset=\Dset$, otherwise. $\DDset$  is equipped with the Skorohod $J_1$-topology and its Borel $\sigma$-field agrees with the product $\sigma$-field $\AA$, see e.g., \cite{MR625374,ferger2015arginf,JJ, MartinE}.
\ED

\BRM  
If $\vk Z \in \DDset$ a.s.\ we retrieve \cite[Lem 2.8]{PH2020} using further the relation between shift-invariant measures $\nu=\nu_Z$ and  {\it shift-generated} $\kk$'s, see \cite{hashorva2021shiftinvariant}. 	
\ERM

For a given CRF $Q$, as shown next $\kkQ$   is shift-generated,  purely dissipative and does not depend on the distribution of the rv $N$. 
\BT \label{conv} Let $\kK$ be given  and let  $\vk Q\in \clasP$. 
If  for some $\TT$-valued rv $N$ as in the Introduction, the RF $Z_N $ defined in \eqref{mm} belongs to $\kK$, then  $Z_N$  belongs to $\kK$ 
for all  rvs $N$ with   pdf $p_N(t)>0,t\in \TT$ being further independent of $\vk Q$. 
Furthermore $\kK$ is shift-generated, $\pk{\SZ \in (0,\IF)}=1$ and  $\vk Q$ is a CRF satisfying further 
\rE{\bqn{\label{supQ} 
	\E*{ \sup_{t\in \TT} \kappa(B^t Q)}\in (0,\IF).
}
}

Conversely, if   $\vk Q\in \clasP$ is a   CRF,  then $\kkQ$ is {\it shift-generated}, 
 \rE{\eqref{supQ} holds} and  
\bqny{\pkT{\SZT \in (0,\IF)}=\pk{\SL \in (0,\IF)}=1 , \quad \forall \ZT \in \kkQ.
}
\ET

The {\it shift-generated} class of the Brown-Resnick RFs  in \Cref{ex:BRModel} satisfies 
$\pk{\kappa(  Z)>0}=1$. The next result explains this in a general framework. 
\BEL Given a {\it shift-generated} $\kk$ the following are equivalent:
\begin{enumerate}[(i)]
\item \label{elker:1} $\kappa( \vk Z )>0$ a.s.;
\item \label{elker:12} For some (and then for all) $\ZT\in \kk$ we have  $\kappa(B^{-t} \ZT )>0$ a.s.\ for all $t\in \TT$;
\item \label{elker:2}  $\kappa( B^{-t}\vk \Theta )>0$ a.s.\ for all $t\in \TT$;
\item \label{elker:3}  $\kappa( B^{-t}\vk Q )>0$ a.s.\ for all $t\in \TT$, provided that 
$\kk=\kkQ$ with   CRF $\vk Q \in \clasP$.
\end{enumerate}
\label{lem:elker}
\EEL

\def\MMM{\mathcal{M}( \vk Y)}	
\def\MMMP{\mathcal{M}_\EI( \vk Y)}

\begin{example} ($m$-truncation of CRFs)
\label{techet}
Let  $\vk Q\in \clasP$ be  a  CRF and fix $m>0$. Setting $\vk Q^{(m)}(t)= \vk Q(t) \ind{ \norm{t} \le m}, t\in \TT$ it follows easily that 
$ \vk Q^{(m)}$ is also a  CRF, provided that $p^{>}_{Q^{(m)}}=1$. Clearly, $\vk Q^{(m)}$ does not satisfy \Cref{lem:elker}, \Cref{elker:3} even when $\vk Q$ satisfies it. Moreover, in general 
$\kkQ$ and $\kkQm$ are different, however by construction they are both purely dissipative.   
\end{example}

\subsection{Constructions of CRFs of purely dissipative $\kk$'s}
Given a purely dissipative    $\kk$, it is of interest to construct CRFs $\vk Q\in \clasP$ such that 
$$\kk= \kkL,$$  
see, e.g.,  \cite{kab2009, dom2016,Htilt,Hrovje,hashorva2021shiftinvariant} for constructions related to max-stable processes and \cite{kulik:soulier:2020,klem,PH2020,MartinE} for new developments and other applications.\\
If the rv $C>0$ satisfies $\E{C}=1$ and $\tQ\in \clasP$ is another CRF, 
then clearly $ Q=C^{1/\alpha} \tQ$ is again a  CRF and by \eqref{cA} 
\bqn{\label{cqq}
\kkQ= \mathcal{C}_{\kappa, N}[\tQ].
	}
 
 If $\tQ^\star \in \clasP$ satisfies the second inequality in \eqref{dense} and 
 $$
 q= 	\pk*{\sup_{t \in \TT} \tQ^\star_\kappa(t)     >0}= \pk*{\sup_{t \in \TT}  \kappa(B^t \tQ^\star)     >0}\in (0,1] , $$
 then  $\tQ^\star$ is not a CRF. For this case, we shall define a CRF $Q$ as follows 
 \bqn{\label{monday} 
 	Q(t) = q^{1/\alpha}\tQ^\star(t) \Bigl \lvert \sup_{s\in \TT} \tQ_\kappa^\star(s)>0, \quad t\in \TT.
	}

As in \cite{kulik:soulier:2020, PH2020} we show next that  CRFs 
can be directly determined by $\vk \Theta,\vk Z $ or $\vk Y$   utilising 
$\mathcal{S}_\EI$ and $\mathfrak{B}_{\EI, \tau}$ (recall \eqref{nowWBack}). 
Set  below   
\bqn{
\MPY=\sup_{t\in \EIE \cap \diad} \kappa^{1/\alpha}( B^{-t} \vk Y)
}
and recall that   
$ \mathcal{P}( \EIE)$ denotes the fundamental parallelogram of the lattice $\EI\subset \TT$ with 
volume $\Delta(\EIE)>0.$

\BT
\label{ffD} 
Let  $\kk$ be   such that $\pk{\SZ <   \IF } =1$. If \eqref{tauR} holds for some 
$\tau\inr$ and $\EIE$ is a  
full rank lattice such that when $\TT=\R^l$ a.e.\  with respect to the measure $\lambda(\cdot)$
\bqn{ \label{doktoret}
 \pkh{ \mathcal{S}_\EI(B^{-t} \vk \Theta) >0 }=1, \quad \forall t\in \mathcal{P}( \EIE) \cap \TT
} 
or $\EIE=\TTT$,  then a  CRF $\vk Q\in \clasP$ such that  $\kk=\kkQ$  can be constructed  
 for all $b\in [ 1,\IF)$ as follows:
\begin{enumerate}[(i)]
\item \label{ffD:1}  $  \vk  Q=   c^{1/\alpha} \vk \Theta ,$ where $1/c=\Delta(\EIE)  \STq$; 
\item \label{ffD:3}    $  \vk Q$ is given by in \eqref{monday}, where  $\tQ^\star=   c^{1/\alpha}
b\vk Y \lvert  \MPY> b,$ with  $ c=    \kappa^{\tau}( Y)/(\Delta(\EIE)\cEE{[\MPY]^\alpha}  \mathfrak{B}_{\EI, \tau}( \vk Y) ) $;
\item \label{ffD:2}   $  \vk Q$ is given  by  \eqref{monday}, where $ \tQ^\star=c^{1/\alpha}Z$ and $c= \frac{ \kappa(\vk Z)}{\Delta(\EIE)   \SZq} $; 
\item \label{ffD:4}  $  \vk Q$ is given by \eqref{monday}, with  $ \tQ^\star=  b \vk Y^{(\tau)} / 
\mathcal{M}_\EIE (\vk Y^{(\tau)})\lvert \mathcal{M}_\EIE(\vk Y^{(\tau)})>b,$ 
where  $\vk Y^{(\tau)}$ is the RF $\vk Y$ under the tilting with respect to $ \kappa^{\tau}( \vk Y)/[\Delta(\EIE)\mathfrak{B}_{\EI, \tau}(\vk Y)]$.
\end{enumerate}
\ET

\def\EJ{\mathcal{Q}}
\def\MTY{ \mathcal{M}(\vk Y)}

It is known from \cite{Hrovje,kulik:soulier:2020, BP, Planic} that 
 CRFs can be constructed by: 
 \begin{enumerate}[(i)]
 	\item utilising shift-involutions acting on $\vk Z$; 
 	\item  positive shift-involutions acting on the {\it spectral tail rf}   $\vk \Theta$;
 	\item anchoring maps applied to the {\it tail rf} $\vk Y$.
 \end{enumerate}
In what follows, we focus on the setting where 
$\EI$ is a full rank lattice on $\TT$ considering positive shift-involutions and anchoring maps denoted by  $\mj_1$ and $\mj_2$, respectively. Further, we shall denote by $\mj_3$ a shift-involution. \\
Note in passing that our definition of anchoring maps is slightly more general than those found in the existing literature.

\BT 
Let $\EI$ be a full rank lattice on $\TT$ \cED{with infinite number of elements} and suppose that $\pk{\SZ< \IF}=1$. If  $\mj_1,\mj_2$ satisfy \eqref{peres2}, $\mj_3$  satisfies \eqref{peres1} and further \eqref{doktoret} holds,   then  
$\kk$ has a  {\it random-shift representation}  with  $\vk Q = c^{1/\alpha} \widetilde{Q}, c>0$   determined  as follows:
\begin{enumerate}[(i)]
\item  \label{invo:1}  $\widetilde{\vk Q}(t)  =  \vk \Theta(t) \lvert (\mj_1(\vk \Theta) =0)$ and   $c=\pkh{\mj_1(\vk \Theta)=0}/ \Delta(\EIE)>0$;
\item  \label{invo:2}  $\widetilde{\vk Q}(t)  =   \frac{\vk Y(t)}{\MPY} 
\lvert (\mj_2(\vk Y) =0, \MPY > b)$ and    $c=b^\alpha 
\pkh{\mj_2(\vk Y)=0, \MPY> b}  /\Delta(\EIE)  , b\in [1,\IF) $;
\item  \label{invo:3}  $\widetilde{\vk Q}(t) =  \vk Z(t) \lvert  (\mj_3(\vk Z)=0)$ and   $c=\pk{\mj_3(\vk Z)=0} /\Delta(\EIE)>0$.
\end{enumerate}
\label{hoff}
\ET 

\BRM \label{fust} 
\begin{enumerate}[(i)]
\item Condition \eqref{doktoret} is fulfilled if $\pk{\kappa( \vk Z)>0}=1$,  since in view of \Cref{lem:elker} this implies $\pkh{\kappa(B^{-t}\vk \Theta)>0}=1$ for all $t\in \TTT$. Hence 
$ \pkh{\mathcal{S}_\EI(B^t \vk \Theta) >0}=1$ follows from \cite[Thm 2.1]{GeoSS}.
\item If $\kk=\kkL$ with $Q$ a CRF, under the assumptions of \Cref{ffD} and \Cref{hoff},  applying \eqref{tcfN} we obtain for all shift-invariant $H\in \Hh_\alpha$
\bqn{\label{stimmt}
	\EHh*{ \frac{H(\vk \Theta)}{\Delta(\EI)\STq}} &=& 
	\EHh*{ \frac{   H(\vk Y)}{\Delta(\EI)  [\MPY]^\alpha} \bigl\lvert J_2(\vk Y)=0}
	=\E*{H(\vk Q) }.
	  \label{stimmt2}
}
\rE{Consequently,  we have 
\bqn{
	\label{simhu}
	\EHh*{ \frac{\max_{t\in \TT} \kappa(B^t \Theta)}{ \Delta(\EI)\STq}} 
	= \frac{\pk{ J_2(\vk Y)=0}}{ \Delta(\EI)}= \E*{\max_{t\in \EI} \kappa(B^t Q) }\in (0,\IF),
}
}
yielding the claims of \cite[Prop 3.9 Eq.\ (23)]{Planic} and \cite[Cor 2.11]{PH2020}.
\item Similar constructions for $d=1$ in both discrete and \cadlag\ case are obtained for slightly less general setting in \cite{debicki2017approximation,Hrovje, kulik:soulier:2020,klem,PH2020}.  
\end{enumerate}	
\ERM

\section{Applications}
\label{sec:cash}
We shall discuss first some properties of the candidate extremal index 
followed by a short investigation on the approximations of purely dissipative $\kk$'s and then continue with an application concerning the Brown-Resnick $\kk$ followed by  a result related to the $m$-approximations.

\subsection{$\EIE$-extremal index}
\label{sec:candi}
Given $\EI$ a lattice on $\TT$ or $\EI=\TT$,  we call  $\vartheta_\EI$ defined by 
\bqn{ 
	\vartheta_\EI
	&=& \EHh*{\frac{1}{\Delta(\EI)\mathfrak{B}_{\EI,0}(\vk Y)}}
}
the $\EI$-extremal index of $\kk$. The finiteness of $\vartheta_{\EI}$ follows from \eqref{leviz}.\\
 We discuss next the case of full rank lattices $\EI$.  In view of \Cref{kabA} 
\bqn{
\vartheta_\EI=0   \iff  \vartheta_\TT=0 \iff \pk{\SZ<\IF}=0,
}
which is also equivalent with one of the events defined in \Cref{kabA} having probability zero.
It follows from the proof of  \Cref{ffD} that  for all $\tau$ satisfying \eqref{tauR}  
\bqn{\label{voice}
\vartheta_\EI
&= & b^{\alpha} 
\EHh*{ \frac{\kappa^\tau(\vk Y)\ind{\mathfrak{B}_{\EI , \tau}(b\vk Y)>0} }
	{   \Delta(\EIE) \mathfrak{B}_{\EI , \tau}(\vk Y) }}< \IF, \quad \forall b\in [1, \IF).
}	
The equality \eqref{voice} for  $b=1,\tau=0$ has been derived  in an unpublished manuscript by the author  under asymptotic restrictions inspired by   \cite{Berman82,Berman92}, where it appears (not explicitly) in relation to the {\it Pickands constants}, see \cite[Thm 10.5.1]{Berman92} and \cite[Thm 1.1]{KSojourn}. We refer to the representation \eqref{voice} of $\vartheta_\EI$ as the {\color{blue} Berman representation}. For $\TT=\Z,b=1,\tau=0$ it  appeared later in \cite{Hrovje,Planic}.  See also \cite{PH2020,hashorva2021shiftinvariant} for a less restrictive framework. 
From the proof of  \Cref{ffD}  it follows that 
\bqn{
\label{lys}
\vartheta_\EI = \EHh*{ \frac{\sup_{t \in \EI} \kappa(B^{-t} \vk \Theta)}{\Delta(  \EI)\STq}}.
}
The representation \eqref{lys}  includes the expression of the extremal index of max-stable stationary RFs  in \cite{Htilt}. As noted in \cite{debicki2017approximation}  that representation is already implied from the seminal papers \cite{Genna04,Genna04c}
and can be thus referred to as the {\color{blue} Samordnitsky  representation}.\\
Other representations for  $ \vartheta_\EI$ are obtained  utilising   \Cref{0LZR}. 
For instance when $\mj_2$ satisfies \eqref{peres2}  
and $\EI$ is a full rank lattice on $\TT$, using further \eqref{0gjelle} 
for all $b\in [1,\IF)$
\bqn{\quad \quad   b^{-\alpha} \Delta(\EIE) \vartheta_\EI  = 
\pkh{\mj_2(\vk Y)=0, \MPY >b}=  	
\pkh*{\sup_{0 \prec  t, t\in \EI} \kappa(B^{-t}\vk Y)\le 1, \MPY>b
}. 
\quad
\label{propprive2}   
} The second expression for $\vartheta_\EI$  above, which follows from the first one taking $\mj_2$ to be the first exceedance map goes back to works of P.\ Albin (case $b=1$ only) and appears as limiting constant (Pickands constant) of supremum of Gaussian and related RFs,  see e.g., \cite{Albin1990,Albin2003}. 

Next, utilising the first exceedance map we obtain from the second expression in \eqref{propprive2}
\bqn{\quad \quad  \quad \Delta(\EIE) \vartheta_\EI =   
	\EHh*{ \sup_{0 \preceq  t, t\in \EI} \kappa(B^{-t}\vk \Theta)  
	-  \sup_{0 \prec  t, t\in \EI} \kappa(B^{-t}\vk \Theta) }=
\E*{ \sup_{0 \preceq  t, t\in \EI} \kappa(B^{-t}\vk Z)  
	-  \sup_{0 \prec  t, t\in \EI} \kappa(B^{-t}\vk Z) }	 \quad 
\label{propprive3}
}
derived for   the  Brown-Resnick max-stable RF in \cite[Corr 6.3]{WS} and initially obtained in \cite{BojanS},  see also \cite{kulik:soulier:2020,HBernulli, BP}. \\
If $\kk=\kkQ$, we have in view of \Cref{ffD} and \eqref{stimmt} the general expression    
\bqn{
	\vartheta_\EI= \E*{ \sup_{t\in \EIE} \kappa(B^{-t}\vk Q)}
}
obtained for $l=1$ in \cite{debicki2017approximation}. 
In view of \eqref{stimmt}, new representations for $\vartheta_\EI$ can be derived by choosing different $Q$'s, see for instance \eqref{simhu}.

In the special case $\kappa( Z)>0$ a.s., and hence \eqref{doktoret} is satisfied, by \Cref{ffD} for $\TT=\R^l$, all full rank lattices $\EI$ and  $\tau$ as in \Cref{tauR} 
\bqn{\label{nastraD1}
\vartheta_\TT=   \frac {b^\alpha}{	\Delta(\EI) }  \EHh*{ \frac{\sup_{t\in \R^l} \kappa(B^{-t}\vk Y) \kappa^\tau (\vk Y) \ind{\MPY >b}} 
	{ \sup_{t\in \EI} \kappa(B^{-t}\vk Y)\mathfrak{B}_{\EI , \tau}( \vk Y) }} , \quad b\ge 1. 
}
\rE{The next example gives an application for the volume of the fundamental  parallelepiped  using the 
first construction in \Cref{ffD}.
\BEX \label{lonay}
  Consider the settings of \Cref{examp:kappa} where we take for simplicity  $d=1$, $Q$ to be non-negative 
   and $\kappa(f)=\abs{f(0)}, \alpha=1$. Suppose further that 
	$\arg \sup_{t\in \TT}  {  Q(t)}=  0$
and set $\kk = \kkL $.  It follows that 
	that $\Theta(t)= Q(t+\mathcal{E}) /{Q(\mathcal{E})},$ 
	with $\mathcal{E}$ having pdf $\Theta$ is the local RF of $\kk$.  Hence for this case for any $\EI$ which is a full rank lattice on $\R^l$ we have  using the first construction in \Cref{ffD}
	\bqny{   
		{  Q(0) } = 	\frac{  \sup_{t \in \TT}  {  Q(t) }}{ 	 	\intT   {  Q(t) }  \lambda(dt) }
		&=& \vartheta_\TT\\
		&=& \frac{1}{ \Delta(\EI)} \E*{  \frac{  \sup_{t\in \TT}  { Q(t- S) }}{
				\sum_{t\in \EIE}   { Q(t- S) } }}
		\\
		&=& \frac{ {   Q(0) }  }{  \Delta(\EI) }  \intT \frac{ {   Q(h)}   } {   \sum_{t\in \EI}    {   Q(t-  h) }   } \lambda(dh)
	}
	implying the following identity for the volume of the fundamental parallellpiped of $\EIE$
	\bqn{\label{inemato}   \Delta(\EI)& =& \intT \frac{    {   Q(h)}   } {    \sum_{t \in \EI}    {  Q(t- h) }   } \lambda(dh) .
	}
	In the special case $Q(t)= e^{-\sum_{i=1}^l\abs{t_i}^2/2} / (\sqrt{2 \pi})^l, t=(t_1\ldot t_l)\inr^l $ and $\EI= \delta \Z^l, \delta >0$ we have 
	$$ \delta^l = \int_{\R^l}  \frac{1}{ \sum_{t\in \delta \Z^l} e^{- \sum_{i=1}^l ( t^2_i/2- h_it_i) }} \lambda(dh), $$
	which has been shown in  \cite{DiekerY} for the case $l=1$. 
	The idea of the above non-asymptotic proof of \eqref{inemato} was kindly communicated by Dima Zaporozhets.
\EEX
}
The relation between $\vartheta_{\TT}$ and $\vartheta_{\EI}$ is first shown for the Brown-Resnick max-stable case in \cite{DiekerY}, see also \cite{ZKE,beman2024}. In view of our results, several other identities relate both constants. In \cite{Hrovje, kulik:soulier:2020} the quantity $\vartheta_{\TT}$ is referred to as the candidate extremal index. As discussed in \cite{BojanS,BojanPhilippe,BP,kulik:soulier:2020,PH2020, drees2022statistical} calculation of the candidate extremal index is closely related to the calculation of extremal index, with few notable exception pointed out  in \cite{wangExt}. Note further that for particular cases, representations of extremal indices are obtained in 
\cite{KW,MR4206416}.
Below we show that $\vartheta_\EI$ is exactly the  extremal index of a corresponding max-stable RF and its calculation can be dealt with within the framework of max-stable RFs. 

\subsection{Max-stable RFs}
\label{sec:maxstable}
Consider a given $\kk$ and the corresponding $\kkA$, which is also {\it shift-generated}. Recall $Z_\kappa(t)=\kappa(B^t Z), t\in \TT$ and hence $Z_\kappa$ is non-negative.  
Let  $ Z^{(i)}_i, i\inn$ be  independent copies of $Z_\kappa$ and define the  max-stable stationary RF $ X_\kappa (t),t\in \TT$   via its de Haan representation  as in \eqref{eq1}. Assume that $\E{Z_\kappa(t)}=1,t\in \TT$, which implies that 
$X_\kappa (t)$ has  a unit $\FRE$ distribution for all $t\in \TT$. Moreover, 
for all $ t_i\in \TT, x_i \in (0,\IF), i\le n$ in view of \cite{MolchanovSPA} 
\bqn{ \label{kode}
	- \ln \pk{ X_\kappa (t_1) \le  x_i, 1 \le i \le n } &=&
	\E*{ \max_{1\le i \le n} \frac {Z_\kappa(t_i)}{x_i^\alpha}     } 
}
and since $Z_\kappa$ has locally bounded sample paths, the law of supremum of $X_\kappa $ on compact intervals is explicitly available, see \eqref{LY} below. In particular, $X_\kappa $ has locally bounded sample paths and  is stationary,  since 
$\kkA$ is shift-generated. \\ 
 It is of interest to derive a {\it  Rosi\'nski representation}  for $X_\kappa $ as in \eqref{eqcR} for a given  CRF $Q$ with corresponding $Q_\kappa$. Such a  representation then 
yields an alternative formula to \eqref{kode}, i.e.,  
\bqn{ \label{kodeQ}
	- \ln \pk{ X_\kappa (t_1) \le  x_i, 1 \le i \le n } &=&
	\E*{ \intT \max_{1\le i \le n} \frac {Q_\kappa(t_i-s)}{x_i^\alpha}   \lambda(ds)  } .
}
\BRM  \begin{enumerate}[(i)]
	\item In view of \eqref{kode} and \eqref{kodeQ},  if $X_\kappa $ has a {\it  Rosi\'nski representation}  with some  CRF $Q_\kappa$, then $\kkA$ has a  {\it random-shift representation} with the same CRF. Moreover, the converse is also true. 
\item If $X_\kappa$ has \cadlag\ sample paths, then we can define the {\it  Rosi\'nski representation}  choosing 
 $\Pi(\cdot)=\sum_{i=1}^\IF \delta_{P_i, T_i, {Q}_{(i),\kappa}}(\cdot)$ on $(0,\IF) \times \mathcal{R}^l \times \DDset$ with mean measure $\lambda_\alpha(\cdot)  \odot c\lambda(\cdot) \odot \mathbb{P}_{{Q}_\kappa^{1/\alpha}} $  (recall $ \lambda_\alpha(dr)= \alpha r^{-\alpha-1} d r$)
 \cED{ and then set 
 \bqn{ X_{\kappa} (t)=  \max_{i\ge 1} P_i B^{T_i}{Q}_{(i),\kappa}^{1/\alpha}(t), \quad t\in 	\TT.  
 } 
}
 For such a choice 
 $$ 1=c\E*{ \int_{\TTT} {Q}_\kappa(t) \lambda(dt) } .    $$
\end{enumerate}
\label{sebPordh}
\ERM
A necessary and sufficient condition for the above mentioned representations is $\pk{\mathcal{S}(Z) < \IF}=1$, which is equivalent with one of the events in \Cref{kabA} holds with probability 1. 
In view of our assumptions  
$$ 
\pk*{\sup_{t\in \TT} Z(t)>0, \mathcal{S}(Z)=0}=0,
$$
hence \cite[Eq.\ 6.5]{Htilt} holds,	 implying that $X_\kappa $ has a Rosi\'nski representation.

The  construction of different $Q_\kappa$'s has been the topic of numerous papers, see e.g.,  \cite{deHaan,kab2009,kabStoev,  dom2016,Htilt, Hrovje,klem,debicki2017approximation, kulik:soulier:2020, PH2020} and the references therein.  Our results imply new constructions when $\tau\not=0, b\in (1,\IF) $ or $\pk{\kappa(\vk Z)>0}=1$.  

\BS  
If $\pk{\mathcal{S}(Z) < \IF}=1$ and $\EI$ is a full rank lattice on $\TT$, then a stochastically continuous  CRF $Q_\kappa$ that defines a {\it  Rosi\'nski representation}  \eqref{eqcR} for $X_\kappa $ can be constructed from \Cref{ffD} or \Cref{hoff}. \label{lemHey} 
\ES
As an application of \Cref{lemHey} we obtain new {\it  Rosi\'nski representation}  of the Brown-Resnick and 
the Brown-L\'evy-Resnick max-stable RFs. For the first case, such representations have been derived in 
\cite{kab2009,debicki2017approximation, Hrovje,kulik:soulier:2020}. The Brown-L\'evy-Resnick max-stable RFs have been studied in \cite{eng2014d,engelke2016levy}.
\BK 
\begin{enumerate} [a)] 
	\item  \label{enmP1}
Let $\kk$ be Brown-Resnick shift-generated. If \eqref{BR2} holds, then $\kk=\kkQ$  and $X_\kappa $ has a {\it  Rosi\'nski representation}  with $Q_\kappa$ determined by \Cref{ffD} or \Cref{hoff};
\item  \label{enmP2} If $\kk$ is a Brown-L\'evy-Resnick class of RFs as in \Cref{exa:lev}, then again   $\kk=\kkQ$ and $X $  has a 
 {\it  Rosi\'nski representation}  with $Q_\kappa$ determined by \Cref{ffD} or \Cref{hoff}. 
\end{enumerate}
\label{korEVG}
\EK
In view of \eqref{kode}, for all $z>0$ and all  full rank lattice $\EIE$ or $\EIE=\TT$ (recall that $X_\kappa $ is also taken to be separable)
\bqn{\label{LY}
-   \ln  \pk*{ \sup_{t \in \cap [0,n]^l \cap  \EIE  } X_\kappa (t) \le  z n^{l/\alpha} }
=  \frac{1}{z^\alpha n^l}   \E*{  \sup_{t \in \cap [0,n]^l \cap \EIE } \kappa(B^{-t}\vk Z )}    =  \frac 1 {z^\alpha} 	\bobo_{  Z}^{\EIE}(n)
}
holds for all $ t_i\in \TT, x_i \in (0,\IF), i\le n, n\inn$.  Consequently, we obtain 
$$\limit{n} \bobo_{Z}^{\EIE}(n)=\bobo_{ Z}^{\EIE}=\vartheta_\EI< \IF.$$
Applying  \cite[Prop 7.2]{hashorva2021shiftinvariant} yields
$\bobo_{Z}^{\EIE}=0$ iff $\pk{\mathcal{S}(Z)=\IF}=1$ and hence  
\bqn{
\bobo_{ Z}^{\EIE} =\pkh{\mathcal{S}(\Theta)< \IF}\limit{n} \bobo_{ Z_*}^{\EIE}(n),   
}	 
with $Z_*$ belonging to the $\alpha$-homogeneous shift-invariant class of RFs generated from the  {\it spectral tail rf}  $ \Theta \lvert \mathcal{S}(\Theta)< \IF$. Consequently, we can assume without loss of generality that $
\pk{\mathcal{S}(Z)< \IF}=1$ and hence 
there exists a CRF $Q$ such that $\kk  = \kkQ$, implying thus 
\bqn{\label{usht}
\bobo_{Z}^{\EIE}= \vartheta_\EI=\E*{ \max_{t\in \EI}  \kappa( B^{-t}Q )}.
}	 
The case $ Z$ has \cadlag\ sample paths  has been considered in  \cite{debicki2017approximation}, see also \cite{kulik:soulier:2020,PH2020,MartinE}. 
\BRM  
	If $Q=c^{1/\alpha}\widetilde{Q}$ with $c>0$ a constant, by \Cref{sebPordh} and \eqref{usht}
\bqn{ \label{remF}
\vartheta_\EI= \frac{ \E{\sup_{t \in \EI} \kappa(B^{-t}\widetilde{Q}) }}
{\E{ \int_{\TT} \kappa( B^{-t}\widetilde{Q}) \lambda(dt)}}\in (0,\IF)
}
and thus we retrieve the claims of \cite[Lem 3, Thm 5]{SBK} when $\EI=\TT$.
\ERM

\begin{example} Let $\alpha=d=1$ and $Z(t)=e^{\overline  W(t)},\overline W(t)=W(t)- Var(W(t))/2, t\inr^l$ be as in \eqref{BR} satisfying  further \eqref{BR2}. Let $\EI=(\delta \Z)^l$ with $\delta >0$ such that $\delta \inn$ if $\TT=\Z^l$ and let $\kappa(f)=\abs{f(0)}$.  Since $ Z(0)=1$ a.s.\ by \Cref{BRownE} $Y(t)= e^{\mathcal  E} Z(t)$ with $\mathcal E$ a unit exponential rv independent of $Z$. For this case we can take $\tau \in [0,\IF)$ and hence for all $b=e^\theta, \theta \ge 0$ and $\Delta(\EI)=\delta^l$, in view of \Cref{ffD}, \Cref{ffD:3} and \Cref{hoff}
\bqn{  
	\vartheta_{\TT}&=& 
	\frac{e^\theta}{\delta^{l}} 
	\E*{  \frac{    \sup_{t\in \TT}  e^{\overline  W(t)} 
			\ind{ \sup_{t\in \EIE} \overline  W(t) +\mathcal  E> \theta} }{ 
			\sup_{t\in \EI} e^{ \overline W(t)}  \sum_{t\in \EI} e^{ \tau \overline W(t)} 
			\ind{ \overline W(t)+ \mathcal  E> 0} 	}}\\
	&			= &\frac{e^\theta}{\delta^l} 			
	\E*{ \frac{ \sup_{t\in \TT} e^{ \overline  W(t)}  } 
		{ \sup_{t\in \EI} e^{ \overline  X(t)} } \ind*{ 
			\sup_{0 \prec  t, t\in \EI} \overline W(t)+\mathcal E\le 0, 
			\sup_{t\in \TT} \overline  W(t)+ \mathcal  E> \theta }}\\
	&=&		e^\theta 
	\E*{ \frac{ \ind{ \sup_{t\in \EI} \overline  W(t)+ \mathcal  E> \theta} } 
		{   
			\intT e^{ \tau \overline W(t)} 
			\ind{ \overline W(t)+ \mathcal  E> 0} 	\lambda(dt)}}\label{cufo1}\\
	&=&	\frac{1}{\delta^l}	 
	\E*{ \frac{  \sup_{t\in \TT} e^{\overline  W(t)} } 
		{   
			\sum_{t\in \EIE} e^{  \overline W(t)}  	}}	
	\label{cufo2}.	
}	
Note that the {\color{blue} Berman representation} \eqref{cufo1} is shown for $l=1$ in \cite{GeoSS} and 
\eqref{cufo2} in \cite{ZKE}.
\label{exampGale}
\end{example}

\subsection{Shift-generated Brown-Resnick $\kk$'s}
The Brown-Resnick max-stable process $X$ that has representer 
$$Z(t)=e^{ W(t)- Var(W(t)/2},\quad t\in \TT= \R,$$ 
with 	$W(t),t\in \TT$ a centered fractional Brownian motion with Hurst parameter $H\in (0,1]$ plays an important role in extreme value theory and statistics. Consider  in this section 
$\kappa(f)=\abs{f(0)}$ 
and recall that we assume the de Haan representation \eqref{eq1} and therefore $X$ has $1$-\FRE marginal df's.\\ 
The case $H=1/2$ is initially studied in \cite{bro1977}, while  $H=1$ has been explored  in \cite{eddy1980distribution,Gale}. Stationarity of $X$ has been established in \cite[Thm 2]{kab2009} for 
 $W$  being Gaussian with stationary increments.  In view of our findings,  the stationarity of $X$ is equivalent with \eqref{tcfN} or \eqref{tcfN:0}. The latter equivalence 
has been shown in \cite{Htilt}, which follows also from previous derivations in \cite{MolchanovBE}. 
\BS Let $\kappa$ be as in \Cref{examp:kappa} and let $\kappa_\star\in \mathcal{H}_\alpha^+$ satisfy \eqref{simhusimhu}. If $\mathcal{C}_{\kappa_\star}[Z]$ is a {\it shift-generated}  Brown-Resnick class or RFs, then it is uniquely defined in terms of the matrix pseudo-cross variogram function $\gamma$. Moreover,
 if  $\mathcal{C}_{\kappa_\star}[Z]$ is purely dissipative, then a  CRF $Q$ that generates this class can be constructed with and its law depends only on $\gamma$ and $\kappa_\star$.
\label{BRprop}
\ES

\subsection{$m$-approximation}
Let below $\kk$ be purely dissipative, i.e., $\pkh{\ST< \IF}=1$ and thus $\kk=\kkQ$ for some  CRF  $\vk Q$.  Define   $\kkQm$ as in \Cref{techet} assuming that $Q^{(m)}$ is a CRF and let  
$$\vk Z^{(m)}_N =  B^N \vk Q^{(m)}/p_N(N)^{1/\alpha}, \quad \vk Z_N =  B^N \vk Q/p_N(N)^{1/\alpha},  $$
with $N$ a $\TT$-valued rv independent of $\vk Q$ with positive pdf $p_N(t)>0, t\in \TT$.
The next result shows that the elements of $\kk$ can be approximated by those of $\kkQm$ as $m\to \IF$, which is in line with the  $m$-approximation developed in \cite{klem,kulik:soulier:2020,BP, PH2020}. 

\BS If $\kappa(0)=0$, then for   all bounded compact sets $K\subset \R^l$ 
\bqn{ \label{MAPROX1} 
	\limit{m}\sup_{n > 0} \frac{1}{n^l} \E*{ \sup_{ t\in n K \cap \TT} \abs{\kappa(B^{-t}\vk Z_N) - \kappa(B^{-t}\vk Z_N^{(m)})}  }
	& =&0.
}
\label{MAPROX}
\ES 

\BRM
Given a full rank lattice $\EI$ on $\TT$ or $\EI=\TT$, if  $\kappa(0)=0$, then 
a direct implication of  
\eqref{MAPROX1} is the following result
\bqn{ 
	\bobo_{\vk Z}^{\EIE} = 	\bobo_{\vk Z_N}^{\EIE}=\limit{m} \bobo_{\vk Z^{(m)}_N}^{\EIE}.
	\label{formB}
}
\ERM

 \section{Proofs}\label{sec:proofs}
 \prooflem{kecat} 
	First note that $X_\kappa$ is stationary, since by \eqref{tcfN:0} we have   
\bqny{
	\E{ Z_\kappa(h)  F(Z_\kappa^{1/\alpha})}=:
	\E{ Z_\kappa(h)  G(Z)}=	\E{ Z_\kappa(0) G(B^h Z)}
	 = 	\E{ Z_\kappa(0) F(B^h Z_\kappa^{1/\alpha})}, 
	\quad \forall	F\in \Hh_{0}, \forall h\in \TT,
} 
which is equivalent with the stationarity of the corresponding max-stable RF $X_\kappa$, see e.g., \cite{Htilt,KumeE}.
Hence, the first claim follows from the characterisation of purely dissipative/conservative max-stable RFs in \cite{dom2016} and is also known from previous works for $\alpha$-stable RFs, see e.g., \cite{Genna04}.  The second claim follows from \eqref{am20} and the assumptions on $\kappa$ using further the equivalence of the norms in $\R^d$. \QED 
\\
\prooflem{springWaltz} The shift-invariance of $\kK$ is an immediate consequence of the stationarity of $Z$. Next, in view of \Cref{kabA} we have that $\pk{\SZ=\IF}=1$ is equivalent with 
$$\pk*{\sum_{t\in \Z^l} \kappa(B^{-t} Z)  = \IF}=1.$$ 
The latter follows from \cite[Cor. 2.1]{ZKE} establishing the claim.
\QED
\\ 
\prooftheo{kabA} 
The claims in \eqref{arm10} and \eqref{am20c} follow from the assumptions \eqref{peres2},\eqref{peres1}   and \Cref{0LZR}.
If $b=1$ and \eqref{lere} holds, the claims  follow from \cite[Lem 6.2, Thm 6.3]{hashorva2021shiftinvariant}. \\
When $b\ge 1$ we have a.s. 
$$\mathfrak{B}_{\EI,\tau}( b\vk Y) \ge \mathfrak{B}_{\EI,\tau}( \vk Y), \quad  \Bigl\{ \lim_{ \norm{t}_* \to \IF, t\in \TT } \kappa(B^{-t}\vk \Theta ) =0   \Bigr \}  \subset \{ \mathfrak{B}_{\EI,\tau}( b\vk Y)< \IF\}
$$ 
and both hold also if $\EI=\TT$. Hence the claims for $b>1$ follow and thus the proof is complete.
\QED  
\\
\proofkorr{kL} In view of \Cref{kecat} and \cite{dom2016}, the pure dissipativity of $\kk$ is equivalent with that of the max-stable RF $X_\kappa$ with 
representer $Z_\kappa^{1/\alpha}$. Since by the assumption $\max_{t\in \TT} Z_\kappa(t)$ is well-defined and a.s.\ finite, then  along the lines of the proof of \cite[Thm 14, p.253]{kab2009}, the pure dissipativity follows if we assume further \eqref{intKab} establishing the proof. 
\QED

\def\SD{S_\diad(Z)}

\prooflem{kor21} 
Taking $F \in \Hh_0$ shift invariant, by \eqref{tcfN:0} if $\EHh{F(\Theta)}=0$ we obtain (set $\SD=\sum_{t\in \diad} p_N(t) \kappa(B^{-t} Z)$) with  $p_N(t)>0, t\in \diad$ summable
\bqny{ 0 &=& \sum_{t\in \diad} p_N(t)\EHh{ F( B^t \Theta)}
	 = \sum_{t\in \diad} p_N(t)\E{ \kappa(Z) F(B^t Z)}\\
&=&\E*{F( Z) \SD }
=\E*{F( \ZT) S_\diad(\ZT ) }  
}
implying that $\E{F(Z)}=0$ since by the assumption $\pk{\sup_{t\in \diad} \kappa(B^{-t} Z) >0} =1$ and thus $$\pk{\SD \in (0,\IF)}=1.$$
Consequently, \Cref{21:a} is equivalent with \Cref{21:b}.\\
\underline{\Cref{21:b} $\implies$ \Cref{21:c}}: If $\E{F(Z)}=0$, then from the above proof $\ETT{F(\ZT)}=0$ for $\ZT \in\kk$. Taking 
    $ \vk Z_N=B^N \vk Q / [p_N(N)]^{1/\alpha} \in \kk$ we have 
    by the shift-invariance and 0-homogeneity of $F$ 
$$
\E{F(\vk Q)}=\E{F(\vk Z_N)}=0
$$
establishing the claim. \\
\underline{\Cref{21:c} $\implies$ \Cref{21:b}}:  If $\E{F(\vk Q)}=0$ holds and $\kk=\kkQ$, then from above  
$\E{F(\vk Z_N)}=0$.\\
When $\kk=\kkQ$, then the stated equivalences imply \eqref{e21}.\\
Conversely, if
$\pkh{ \ST \in (0,\IF)}=1$ we have that $\vk Q= c^{1/\alpha}\vk \Theta, c=1/\ST$
is well-defined.   For all $F\in \wH_0$, using that $\pk{\kappa(\Theta)=1}=1$ and applying 
\eqref{eqDo20}
 we have for 
$ Z_N(t) = 
	B^N \vk Q(t) [p_N(N)]^{-1/\alpha} , t\in \TT$
\bqny{  \E{ \kappa(B^{-h} Z_N )} &=&
	\intT \EHh*{  \frac{\kappa(B^{-h+t} \Theta)}{\ST} F(B^{t}\Theta)}\lambda(dt)\\
	&=&  \intT \EHh*{  \kappa(B^{-h+t} \Theta) \frac{\kappa(\Theta)}{\ST }F(B^t\Theta)}\lambda(dt)
	\\
		&=&  \intT \EHh*{   \frac{\kappa(B^{t-h} \Theta)}{\ST }F(B^h\Theta)}\lambda(dt)
		\\
&=&  \EHh*{   F(B^h \Theta) 
	\intT  \frac{\kappa(B^{s} \Theta)}{\ST }\lambda(ds)}\\
&=&  \EHh{F(B^h \Theta)}= \E{ \kappa(Z) F(B^h Z)}, \quad \forall F\in \wH_0
}
implying \eqref{tcfN:0} and hence $Z_N$ belongs to $\kk$ establishing the proof.
\QED

\prooftheo{conv}  
If for some $\TT$-valued rv $N$ with pdf $p_N(t)>0,t\in \TT$ independent of $\vk Q$  we have 
$Z=\vk Z_N=p_N(N)^{-1/\alpha} B^N \vk Q \in \kk$, then  applying \eqref{tcfN:0} and the Tonelli Theorem,  
we obtain  for all $k\in \TT$
\bqny{ \EHh{F(B^h \vk \Theta)}&=& \E{ \kappa(B^{-h}\vk Z) F(\vk Z)}\\
&=& \E*{ \intT \kappa( B^{t-h} \vk Q)  F(B^t \vk Q) \lambda(dt)  } \\
&=& \E{ \kappa(  B^{-h}B^k\vk Z_{N^*})F(B^k\vk Z_{N^*})}, \quad \forall F\in \Hh_{0},
}	
where the first equality follows since by the assumption $\vk Z_N \in \kK$ and $N^*$ is another $\TT$-valued rv with positive pdf $p_{N^*}$ being further independent of $\vk Q$. This shows that 
$$B^k\vk Z_{N^*} \in \kK$$
 independent of the choice of $N^*$. In particular $B^k Z \in \kK$ for all $k\in \TT$ and thus $\kK$ is shift-generated. \\
Next, by \eqref{perfect} and the shift-invariance of the Lebesgue measure  
$$ 
1=\E{ \kappa( \vk Z) } =
\E*{  \intT   \kappa( B^{-t}  \vk Q)\lambda(dt) }  = \E{ \SL} \implies \pk{\SL < \IF}=1$$
 and thus from \eqref{am0} and \Cref{kor21} $\pk{\SZ \in (0, \IF)} =1$. 
Further we have  
$$
1=\pk*{ \sup_{t\in \TT} \kappa(B^{-t} Z)>0} 
=\pk*{ \sup_{t\in \TT} \kappa(B^{-t}Q)>0}=:p_{Q_\kappa}^>.
$$
Suppose for simplicity in the rest of the proof that $l=1$ and let $M$ be a positive integer.   
By  \eqref{dense} and \eqref{mm},  the Tonelli Theorem implies 
\bqny{  
\IF > \E*{ \sup_{t\in [0,2M+1] \cap \TTT} \kappa(B^{-t} \vk Z)} & =& 
\E*{ \intT  \sup_{t\in [0,2M+1] \cap \TTT} \kappa(B^{-t-x} Q)  \lambda(dx)  }\\
&=&  \sum_{i\in \Z} \int_{i}^{i+1} \E*{  \sup_{t\in [0,2M+1] \cap \TTT }\kappa(B^{-t-x} Q)    } \lambda(dx)\\
&\ge &  \sum_{i\in \Z} \int_{i}^{i+1} \lambda(dx) \E*{  \sup_{s\in [i+1+M, i+1+ 2M] \cap \TTT} \kappa(B^{-s} Q)   }\\
&= &  \sum_{j\in \Z}  \E*{  \sup_{s\in [j, j+ M] \cap \TTT} \kappa(B^{-s} Q)    }\\
&\ge & \E*{  \sup_{s\in \TTT} \kappa(B^{-s} Q)    },
}	  
hence \eqref{supQ} holds. Furthermore, by the above derivations and the shift-invariance of the measure $\lambda(\cdot)$
\bqny{ \IF &>& \E*{ \intT  \sup_{t\in [0,2M] \cap \TTT} \kappa(B^{-t-x} Q)  \lambda(dx)  }
	= \intT \E*{    \sup_{t\in [-M,M] \cap \TTT} \kappa(B^{-t+x} Q) } \lambda(dx) \\
&=&\intT \E*{    \sup_{t\in [-M,M] \cap \TTT} Q_\kappa(t-x) } \lambda(dx)
}
 the second condition in \eqref{dallandysh} follows, implying that $ Q$ is a  CRF.

To prove the converse, assume that  
$	 \vk Z  =\vk Z_N$ is given as above with  $\vk Q$  a  CRF. From the above derivation, also $B^k Z_{N}$  belongs to $\kkQ$ for all $k\in \TT$ and thus $\kkL$ is shift-generated. 
Since by the assumption $p_{Q_\kappa}^>=1$, the independence of both $N, Q$  and the fact that $Q\in \clasP$ imply 
$$
\pk*{\sup_{t\in \TT} \kappa(B^{-t} Z )>0}=\pk*{\sup_{t\in \diad} \kappa(B^{-t} Z )>0}=1.
$$ 
Moreover, in view of \eqref{boll} and the independence of $N$ and $Q$ we have 
\bqny{  
\E*{ \sup_{t\in [-c,c] \cap \TTT} \kappa(B^{-t} Z)} & =& 
\E*{ \intT  \sup_{t\in [-c,c]\cap \TTT} \kappa(B^{-t-x} Q)  \lambda(dx)  } < \IF, \quad \forall c>0.
}	  
By the definition we have further 
\bqn{ \label{wieder} 
\E{\kappa(B^{-h} Z)}=\E{\SL} \in (0, \IF) , \quad \forall h\in \TT
}	
implying thus \eqref{dense} and hence $\kkL=\kK$ is a well-defined $\alpha$-homogeneous class of RFs from $\clasP$. The rest of the proof follows from the arguments given in the the first part above. 
\QED 

\prooflem{lem:elker} The first three equivalences follow from \cite[Lem 9.7]{hashorva2021shiftinvariant}. Assume next that $\kk=\kkQ$ with $\vk Q\in \clasP$ a  CRF. 
If \Cref{elker:3} holds and thus $\kappa(B^{-t} Q)>0$ a.s.\ for all $t\in \TT$, then since $\vk Z_N$ defined by \eqref{mm} belongs to $\kk$ by the $\alpha$-homogeneity of $\kappa$  we have a.s.\ 
$$
\kappa(B^{-t}\vk Z_N)= \frac{1}{p_N(N)} \kappa(B^{N-t} Q)>0, \quad \forall t\in \TT
$$ 
implying \Cref{elker:12}. \\
 Next, if \Cref{elker:2} is valid, then  taking $\vk Q= c^{1/\alpha}\vk \Theta, c=1/\ST$  \Cref{elker:3} follows. Note that $Q$ is a valid CRF, which is consequence of \Cref{kabA} and \eqref{tcfN:0}. See also the first claim in \Cref{ffD}. 
\QED
 
\prooftheo{ffD}   
It suffices to show that  for $\vk Z_N=  p_N(N)^{-1/\alpha}  B^N \vk Q$ and for all $F\in \Hh_0, h\in \TTT$
\bqn{	\label{eHG}
\E{ \kappa(B^{-h} \vk Z_{N})  F(\vk Z_{N})}=	\E{ G_h ( \vk Q)  } = \EHh{ F(B^h \vk \Theta)},
}
with
\bqn{ \label{Gh}
G_h( f)=  \intT  \kappa( B^{y-h} f)  F( B^y  f ) \lambda(dy), \quad f\in \Dset.
}	
Note in passing that by definition and the $0$-homogeneity of $F$ we have 
	$$\EHh{ F(B^h \vk \Theta)}=\E{\kappa(B^{-h} \vk Z) F(\vk Z)}.$$
\\
\underline{Proof of \Cref{ffD:2}, \Cref{ffD:1}}: It is enough to show the proof of \Cref{ffD:1}.  
Since $\pkh{\ST \in (0,\IF) }=1$, which by \Cref{kabA} is equivalent with $\pkh{\STq \in (0,\IF)}=1$, then the RF   
$\vk Q(t)= c^{1/\alpha}  \vk \Theta(t),t\in\TT$ with $c=1/(\Delta(\EI) \STq)$  belongs to $\clasP$. For $F\in \Hh_0$ and $h\in \TT$
we have 
\bqn{ 
\EHh{ G_h ( \vk Q)  }  	&=& \EHh*{   \frac{ G_h(\vk \Theta) }{ \Delta(\EIE)\STq}  }.
\label{dsh}
}	
Note that when $\EI=\TTT$ we set $\Delta(\EIE)=1$. For  this case  applying 
\eqref{eqDo20},   
for all   $\Gamma\in  \Hh_0,  h\in \TT$  we obtain 
\bqn{ 
\EHh*{   \frac{ G_h(\vk \Theta) }{ \ST}  }	&=&\EHh{ F(B^h \vk \Theta)}
\label{dsh2}
}	
and hence \eqref{eHG} follows. Consider therefore next  the case $\EI$ is a full rank lattice of $\TT=\R^l$ and thus 
\bqn{\label{span} 
span(\EIE)=span(\Z^l)=\R^l,
}
where $span(E)$ denotes the smallest linear subspace of $\R^l$ containing $E \subset \R^l$.  It is well-known (see e.g.,  \cite[Lem 10.7]{MR3887637}) that we can tilt $\TT$ by the fundamental domain of the full rank lattice $\EI$ on $\TT$. 
We take as fundamental domain the fundamental parallelpiped  $\mathcal{P}( \EIE)= \{ A x, x\in [0,1)^l\} $, where $A$ is a $l\times l$ base matrix which is non-singular since $\EIE$ is a full rank lattice; recall  $\EIE= \{ A x, x\in \Z^l\}$ and 
$$Vol(\mathcal{P}(\EIE))= \Delta(\EIE) = 
det(A)>0.$$
Consider $\TT=\R^l$ which is spanned by $\EIE$. Hence we have the  tiling  of $\TT$ as 
$$\TT=span(\EI)=\cup_{t\in \EIE}\{ t+ \mathcal{P}( \EIE) \},$$
where  
$t+ \mathcal{P}( \EIE)$ and $s+ \mathcal{P}( \EIE)$ are disjoint for $t\not=s \in \R^l$. 
By \Cref{kabA} and \eqref{doktoret}   
\bqny{ \mathcal{A}(\vk \Theta)= \{  \ST \in (0, \IF),  \STq \in (0, \IF) \} = \{   \ST\in (0, \IF) \}=\{   \STq \in (0, \IF)\} ,\\
 \mathcal{A}(B^{s+t}\vk \Theta)= \{  \ST \in (0, \IF),  \mathcal{S}_\EI(B^{s+t} \vk \Theta) \in (0, \IF) \} =\mathcal{A}(B^{t}\vk \Theta)= \{  \ST \in (0, \IF),  \mathcal{S}_\EI(B^t \vk \Theta) > 0) \}
}
for all $ s \in \EI, t\in \mathcal{P}( \EIE) \cap \TT$. Further,  \eqref{doktoret} yields a.s.
$$ 
\mathcal{A}(B^{s+t}\vk \Theta)=\mathcal{A}(B^{t}\vk \Theta)=\mathcal{A}(\vk \Theta).
$$
Write hereafter $\E{ A;B}$ instead of $\E{A \ind{B}}$ and take $F \in \Hh_0.$ \\
 Using \eqref{eqDo20} for the derivation of the fourth line below, the Tonelli Theorem and the shift-invariance of $\lambda(\cdot)$ 
(recall that $ \kappa( \vk \Theta) =1$ a.s.\ and we interpret $0:0$ as 0)
\bqn{ 
\lefteqn{\EHh*{ \frac{ G_h(\vk \Theta)}{ \STq}}}\notag\\
&=& \EHh*{ \frac{ G_h(\vk \Theta)}{ \STq}; \mathcal{A}(\vk \Theta)} \notag\\
&=&\EHh*{  \frac{\ST}{ \ST}\frac{ G_h(\vk \Theta)}{ \STq}; \mathcal{A}(\vk  \Theta)} \notag \\
&=&   \int_{ t-s \in  \mathcal{P}(\EIE) } \sum_{s\in \EI}   \EHh*{ \frac{ \kappa( B^{-t}\vk \Theta)  }{ \STq}\frac{ G_h( \vk \Theta)}{ \ST} ;  \mathcal{A}(\vk \Theta)} \lambda(dt)  \notag \\
&=&   \int_{ t \in  \mathcal{P}(\EIE) } \sum_{s\in \EI}  \EHh*{ \kappa( B^{-s-t} \vk \Theta)  \frac{ G_h( \vk \Theta)}{ \ST} \frac{  \kappa( \vk \Theta) }{ \STq}; \mathcal{A}(\vk \Theta)} \lambda(dt)  \notag  \\
&=&   \int_{ t \in  \mathcal{P}(\EIE) } \sum_{s\in \EI}  \intT 
\EHh*{ \kappa(B^{-s-t} \vk \Theta)  \frac{ \kappa(B^{y-h}  \vk \Theta)  F(B^y \vk \Theta)}{ \ST} \frac{  \kappa( \vk \Theta)  }{ \STq}; \mathcal{A}(\vk \Theta)} \lambda(dy)\lambda(dt)  \notag  \\
&=&   \int_{ t \in  \mathcal{P}(\EIE) } \sum_{s\in \EI}  \intT 
\EHh*{  \frac{ \kappa(B^{s+t+y-h}  \vk \Theta)  F(B^{s+t+y} \vk \Theta)}{ \ST} \frac{  \kappa(B^{s+t}\vk \Theta) }{ \mathcal{S}_\EIE (B^{s+t} \vk \Theta )}; \mathcal{A}(B^{s+t}\vk \Theta)} \lambda(dy)\lambda(dt)  \notag  \\
&=&  \int_{ t \in  \mathcal{P}(\EIE) }  \sum_{s\in \EI}
 \EHh*{    \frac{ G_h( \vk \Theta)}{ \ST} \frac{ \kappa(B^{s+t}\vk \Theta) }{ \mathcal{S}_\EIE (B^{s+t} \vk \Theta )} ; \mathcal{A}(B^{s+t}\vk \Theta)   } \lambda(dt)  \notag \\
&=&  \int_{ t \in  \mathcal{P}(\EIE) }  \EHh*{    \frac{ G_h( \vk \Theta)}{ \ST} \frac{ \sum_{s\in \EI} \kappa(B^{s+t}\vk \Theta) }{ \mathcal{S}_\EIE (B^{t} \vk \Theta )} ; \mathcal{A}(B^{t}\vk \Theta)    } \lambda(dt)  \notag \\
&=& \int_{ t \in  \mathcal{P}(\EIE) }   \EHh*{    \frac{ G_h( \vk \Theta)}{ \ST} \frac{\mathcal{S}_\EIE (B^{t} \vk \Theta )}{ \mathcal{S}_\EIE (B^{t} \vk \Theta )} ; \mathcal{A}(B^{t}\vk \Theta)  } \lambda(dt)  \notag \\
&=& \int_{ t \in  \mathcal{P}(\EIE) }   \EHh*{    \frac{ G_h( \vk \Theta)}{ \ST}; \mathcal{A}(\vk \Theta)     } \lambda(dt)  \notag \\
&=& \Delta(\EIE)\EHh*{       \frac{ G_h( \vk \Theta)}{ \ST}} \notag\\
&=& \Delta(\EIE)\EHh{ F(B^h \vk \Theta)},
\label{LHope}
}
where the last equality  is implied by \eqref{dsh}. 
\rE{Hence \eqref{eHG} is valid and thus the claim follows. }\\
Next, take $\TT=\Z^l$ and  define the additive quotient   group $\TT/ \EIE=\{ x+ \EIE, x\in \TT\}$. In view of \eqref{span}
we have that the order $m$ of the quotient group is given by (see \cite{Dadush})
$$m=\abs{ \TT/\EIE}= \abs{ \Z^l \cap \mathcal{P}(\EIE)}= \Delta(\EIE)/\Delta(\TT)=\Delta(\EIE).$$
This shows that we have a tiling  of $\TT$ by $\TT \cap \mathcal{P}(\EIE$). Hence repeating the above calculations by substituting integration with summation  establishes the claim.  \\
\underline{Proof of \Cref{ffD:3}}: 
Recall that in our notation $\MPY= \max_{t\in \EI} \kappa^{1/\alpha}(B^{-t}\vk Y)  $. 
  Since $\pkh{\ST < \IF}=1$, then \Cref{kabA} and \nelem{molla} show that  
$\vk Q$ is well-defined and belongs to $\clasP$. For this choice of $\vk Q$ and 
 $G_h$ as in \eqref{Gh} for all $b\ge 1$ we have  
\bqn{ 
\EHh{ G_h(\vk Q)} 	&=& b^\alpha
\EHh*{   \frac{ G_h(\vk Y)  \kappa^\tau(Y) \ind{\MPY> b}}{ \Delta(\EIE) [\MPY]^\alpha \mathfrak{B}_{\EI, \tau}( \vk Y)}  }.
\label{LYL}
}	
By \Cref{kabA} and  $\pkh{\MPY>1}=1$ we have modulo null sets for all $b\ge 1$ (recall that $\kappa(Y)>1$ a.s.)
\bqn{ \label{Ay} 
	\mathcal{A}(\vk Y )&=& 
\{ \mathcal{S}_\EIE(\vk Y) \in (0,\IF)\}    
=\{\mathcal{S}_\EIE(\vk Y) \in (0,\IF),  \mathfrak{B}_{\EI, \tau}(\vk Y) \in (0,\IF)\}=\mathcal{A}(b\vk Y).
}
Moreover, in view of \cite{GeoSS} we have $\{ \MPY> b\} \subset \mathfrak{B}_{\EI, \tau}( \vk Y/b)  $ for all $b\ge 1$. Further 
$
\mathfrak{B}_{\EI, \tau}( \vk Y/b) \le \mathfrak{B}_{\EI, \tau}( \vk Y)< \IF$ a.s.\ and thus by  
the Tonelli Theorem, the shift-invariance of $\lambda(dt)$ and \eqref{tYY}  applied to obtain the second last equality, we obtain 
\bqn{ \lefteqn{
\EHh*{ \frac{ G_h(\vk Y ) \kappa^\tau(Y)  \ind{\MPY> b  }}{ [\MPY]^\alpha \mathfrak{B}_{\EI, \tau}( \vk Y)  } } } \notag \\ 
&=&   \EHh*{ \frac{ G_h(\vk Y ) \kappa^\tau(Y)  \ind{\MPY> b  }}
{ [\MPY]^\alpha\mathfrak{B}_{\EI, \tau}( \vk Y)};\mathcal{A}(\vk Y) }
\notag\\  
&=&    \EHh*{ \frac{ G_h(\vk Y ) \kappa^\tau(Y)  \ind{\MPY> b  } \mathfrak{B}_{\EI, \tau}( \vk Y /b)}
{[\MPY]^\alpha \mathfrak{B}_{\EI, \tau}(\vk Y)\mathfrak{B}_{\EI, \tau}( \vk Y /b)} ;\mathcal{A}(\vk Y) }
\notag\\
&=&    \EHh*{ \intDD \frac{ G_h(\vk Y ) \kappa^\tau(Y)  \ind{\MPY> b  } 
\ind{\kappa(B^{-t}Y/b) >1  } }{[\MPY]^\alpha \mathfrak{B}_{\EI, \tau}( \vk Y)  
\mathfrak{B}_{\EI, \tau}( \vk Y /b)};\mathcal{A}(\vk Y)}\lambda(dt)\notag\\
&=&   \EHh*{ \intDD \frac{ G_h(\vk Y ) \kappa^\tau(Y)  \kappa^\tau( B^{-t}\vk Y/b)
\ind{\kappa( B^{-t} Y /b )>1 } }{[\MPY]^\alpha \mathfrak{B}_{\EI, \tau}( \vk Y) 
\mathfrak{B}_{\EI, \tau}( \vk Y/b)}; \mathcal{A}(\vk Y) }\lambda(dt)\notag\\
&=&   \intDD \intT \EHh*{  \frac{ \kappa(  B^{y-h} \vk Y)F(B^y \vk Y)   \kappa^\tau(B^{-t}  Y/b) \kappa^\tau(Y)
\ind{\kappa( B^{-t} Y/b)>1  } }{[\MPY]^\alpha \mathfrak{B}_{\EI, \tau}( \vk Y) 
\mathfrak{B}_{\EI, \tau}( \vk Y/b)}; \mathcal{A}(\vk Y) } \lambda(dy) \lambda(dt)\notag\\
&=&   b^{-\alpha}\EHh*{ \frac{G_h(\vk Y )  \kappa^\tau(Y)}{[\MPY]^\alpha\mathfrak{B}_{\EI, \tau}( \vk Y) }   \intDD 
	\frac{ \kappa^\tau (b B^t Y)
\ind{\kappa( b B^{t} Y)> 1  } }{	\mathfrak{B}_{\EI, \tau}( b\vk Y)  
}\lambda(dt);\mathcal{A}(b \vk Y)} \notag\\
&=&   b^{-\alpha}\EHh*{ \frac{G_h(\vk Y )  \kappa^\tau(Y)}{[\MPY]^\alpha  \mathfrak{B}_{\EI, \tau}( \vk Y)} }
\label{deshtimos2}
}
holds for all $F\in \Hh_0$ shift-invariant with respect to $\EIE$. For the second last equality we also used \eqref{Ay}.

From \eqref{Ay} and as shown in \cite{hashorva2021shiftinvariant} since $\pk{\mathcal{A}(\vk Y ) }=1$, then a.s. 
$$ \mathcal{A}(\vk Y )=\mathcal{A}(\vk Y/z )=
\mathcal{A}(\vk Y) \cap \{ \mathfrak{B}_{\EI, \tau} ( \vk Y / (z\MPY)) \in (0,\IF)\} =\mathcal{E}_z
$$
for each fixed $z\in (0,1)$ up to a set  with Lebesgue measure zero.  Using the Tonelli Theorem and \eqref{tYY} for the derivation of the third equality below, we have  
(recall   $\lambda_\alpha(dz)=\alpha z^{- \alpha -1} dz$ and set $z_Y= z \MPY$)  
\bqny{ 
\lefteqn{  \EHh{ G_h(\vk Q)}  } \\
&=&   \int_{\EIE} \EHh*{ \frac{ G_h(\vk Y) \kappa^\tau(Y)}{\mathcal{S}_\EIE(\vk Y) [\MPY]^\alpha } \frac{ \kappa(B^{-h} Y) }{  \mathfrak{B}_{\EI, \tau}( \vk Y) };\mathcal{A}(\vk Y ) }\lambda(dh) \\
&=&   \int_{\EIE}  \int_{0}^\IF   
 \EHh*{  \frac{ G_h(\vk Y) \kappa^\tau(Y)}{ \mathcal{S}_\EIE(\vk Y)   [\MPY]^\alpha} \frac{ \ind{   \kappa(zB^{-h} Y) > 1}}{ \mathfrak{B}_{\EI, \tau}( \vk Y )};\mathcal{A}(\vk Y ) }   
\drz  \lambda(dh) \\
&=& \int_{\EIE}  \int_{0}^\IF    
\EHh*{ \frac{ G_h(\vk Y) \kappa^\tau(B^h Y/z) }{ \mathcal{S}_\EIE(\vk Y)    [\MPY]^\alpha } 
\frac{	\ind{    \kappa(B^h \vk  Y/z) > 1}}{ \mathfrak{B}_{\EIE,\tau}( \vk Y /z) };\mathcal{A}(\vk Y/z ) } 	\drzz \lambda(dh) \\
&=& \int_{\EIE}  \int_0^{1}    
\EHh*{ \frac{ G_h(\vk Y) \kappa^\tau(B^h  Y/z_Y) }{ \mathcal{S}_\EIE(\vk Y)    } 
\frac{	\ind{    \kappa(B^h Y/ z_Y)>1  }}{ \mathfrak{B}_{\EIE,\tau}( \vk Y / z_Y )   }; \mathcal{E}_z  } 	\drzz  \lambda(dh) \\
&=& \int_0^{1} \EHh*{  \frac{ G_h(\vk \Theta)}{ \STq  }   \frac{ \int_{\EIE}  \kappa^\tau ( B^h  Y /z_Y ) \ind{  
			\kappa( B^h \vk Y/z_Y) >   1  }  \lambda(dh) }{\mathfrak{B}_{\EIE ,\tau} (    \vk Y /z_Y) }; \mathcal{E}_z }  
\drzz  \\
&=&    \int_0^1  \EHh*{ \frac{ G_h(\vk \Theta) } { \STq} ; \mathcal{A}(Y)}
\drzz  \\
&=&   \EHh*{ \frac{ G_h(\vk \Theta) } { \STq} },
}	
where 1 in the upper bound of the integral in the last fourth  line above is  justified by the fact that a.s.
$$
\ind{\kappa(B^{t} \vk Y /(s\MME))>1 }=0, \quad  \forall s>1,
$$
 whereas the $\alpha$-homogeneity of $\kappa$ was used for the derivation of the third equality above.  
Hence in view of \eqref{LHope} we establish \eqref{eHG}.  \\
\underline{Proof of \Cref{ffD:4}}: We have that 
$$\vartheta_\EI =\EHh{   \kappa^\tau(Y) /(\Delta(\EI)\mathfrak{B}_{\EI, \tau}(\vk Y))} \in (0,\IF),$$
 which is clear for $\EI$ being a lattice on $\R^l$, while for $\EI=\TT$ it is a consequence of \eqref{leviz}. Taking $\vk Y^{(\tau)}$ having the same law as   $\vk Y$ under  
$$
\mathbb{P}^*\{A\}= \vartheta_\EI^{-1}\EHh{ \kappa^ \tau( Y)  /(\Delta(\EI) \mathfrak{B}_{\EI, \tau}(\vk Y ) ) \ind{  A}}, \quad A \in \mathscr{F}
$$
the claim follows from the calculations in the previous case. 
\QED 

\def\EI{\EIE}
\def\JJxi{\mathbb{B}_\EIE(\vk Y/x)} 
\def\MPYn{ \mathcal{M}_{\EIE_n}(\vk Y)}
\def\JJqn{\mathfrak{B}_{\EIE_n}(\vk Y)}
\def\TT{\EIE}
\def\TT{\mathcal{T}}	

\prooftheo{hoff} First note that by the assumption $\pkh{\ST < \IF}=1$. Hence \Cref{kabA} implies 
$$
\pkh{\STq \in (0,\IF)}= \pkh{ \JJq \in (0,\IF)}=1
$$
and moreover, since $\vartheta_\EI \in (0,\IF)$ we have from 
\eqref{0e11A}, with $F(f)= \max_{t\in \EI} \kappa(B^t f)$
\bqn{ \label{0e11A:B} 
\quad\quad 0< 	\EHh*{  \frac{F(\vk \Theta  )}{\STq} } =
	\EHh*{  F(\vk \Theta)  ; \mj_1(\vk \Theta) =0  }
	=	
	\EHh*{  F(\overline{\vk \Theta})  ; \mj_2(\vk Y) =0  } =		
	\E*{  F(\vk Z) ;  \mj_3 (\vk Z)=0  } < \IF.
}
The positivity of the expressions in \eqref{0e11A:B} is a consequence of \eqref{simhu}.
Using further \eqref{0gjelle:2} it follows that the constant $c$ is positive in both three cases treated below.\\
\underline{Proof of  \Cref{invo:1}}: With the notation of \Cref{conv} 
we have 
\bqn{ \label{0e11A} 
	\EHh*{  \frac{F(\vk \Theta  )}{\STq} } =
	\EHh*{  F(\vk \Theta)  ; \mj_1(\vk \Theta) =0  }
	=	
	\EHh*{  F(\overline{\vk \Theta})  ; \mj_2(\vk Y) =0  } =		
	\E*{  F(\vk Z) ;  \mj_3 (\vk Z)=0  } 
}
and using \eqref{0e11A} we obtain 
$$ 
\EHh{ G_h ( \vk Q)  }= \frac{1}{\Delta(\EI) } \E{ G_h ( \vk \Theta); \mj_1(\vk \Theta) =0  }= 
\EHh*{  \frac{G_h(\vk \Theta  )}{\Delta(\EI) \STq} }.  
$$
Hence \eqref{LHope} implies  \eqref{eHG}   establishing   the claim.\\ 
\underline{Proof of  \Cref{invo:2}}: Applying  \eqref{0gjelle} and then \eqref{deshtimos2}    
we obtain  (recall $\MPY= \sup_{t\in \EIE} \kappa^{1/\alpha}(B^{-t}  Y)$)  
\bqny{ 
\EHh{ G_h ( \vk Q)  }
&=& b^\alpha 
\frac{1}{\Delta(\EI) } \EHh*{ \frac{ G_h ( \vk Y) \ind{ \MPY> b}}{ [\MPY]^\alpha}; \mj_2(\vk Y) =0  }\\
&=& b^\alpha 
\frac{1}{\Delta(\EI) } \EHh*{ \frac{ G_h ( \vk Y) \ind{ \MPY> b}}{ [\MPY]^\alpha \mathfrak{B}_{\EI,0}(\vk Y)}  }\\
&=&\frac{1}{\Delta(\EI)} \EHh*{  \frac{G_h(\vk Y  )}
{ 
[\MPY]^\alpha \mathfrak{B}_{\EI,0}(\vk Y) } } \\
&=& \EHh*{  \frac{G_h(\vk \Theta  )}{\Delta(\EI) \STq} },
}
where the last equality is shown in the proof of \Cref{conv}. 
Hence again  \eqref{eHG} is satisfied  and thus  the claim follows.\\
\underline{Proof of  \Cref{invo:3}}: The proof is established by applying \eqref{0e11A}.    
\QED

\def\JJxi{\mathbb{B}_\EIE(\vk Y/x)} 
\def\MPYn{ \mathcal{M}_{\EIE_n}(\vk Y)}
\def\JJqn{\mathfrak{B}_{\EIE_n}(\vk Y)}

\COM{\prooflem{gjakitij}
\Cref{conv} implies $\SQ\in (0,\IF)$ a.s.\ and hence for $Z$ such that 
$\kk= \kkL$
the shift-invariance of $F \in \wH_\alpha$ yields
\bqny{
	\E*{\kappa(\vk Z) \frac{F(\vk Z)}{\SZ}}=\E*{ F(Q) \int_{\TT}   \frac{\kappa(B^h\vk Q)}{\mathcal{S}(Q)} \lambda(dh)}
	=\E*{F(\vk Q) }
}
establishing the claim. \\
We prove next the converse claim. If $Z \in \kk= \kkL$, then \eqref{anabechova} yields for all $F\in \TwH$ 
\bqny{
	\EHh{F(\Theta)}= \E{ \kappa(Z) F(Z/\kappa^{1/\alpha}(Z ) )}
	=\E*{  H(Q)} ,   \quad H(f)=\int_{\TT} \kappa(B^h f) F(B^h f /\kappa^{1/\alpha}(B^h f) ) \lambda(dh).
}
Since by the translation invariance of the Lebesgue measure, the map $H$ defined above is shift-invariant and moreover $H\in \Hh_\alpha$ because $\kappa\in \Hh_\alpha$, 
then the assumption that $\E{H(Q)}=\E{H(\tQ)}$ for all $H \in  \Hh_\alpha$ implies 
$$ 
\E{ \kappa(Z) F(Z/\kappa^{1/\alpha} (Z )) } = \ETT{ \kappa(\ZT) F(\ZT/\kappa^{1/\alpha} (\ZT ))}, 
\quad \forall \ZT \in \mathcal{C}_{\kappa,N}[  \tQ ].
$$
Consequently,  the claim follows from \Cref{remMaps} and \Cref{remMapsP}. \QED 
}

\proofprop{lemHey} 
Since $\kk$ is shift-generated, by \cite[Lem 7.1]{hashorva2021shiftinvariant}  $  X $ is stationary. The assumption \eqref{boll} and \eqref{kode} implies that it has locally bounded sample paths and it is stochastically continuous. 
As shown in  \cite{dom2016}, in view of  \cite[Lem 2]{deHaan} $ X $ has a representer $Z^*$ which is stochastically continuous. In view of \cite[Thm 2.6]{Htilt} and \eqref{tcfN:0} it follows that $Z^* \in \kk$. As in the proof of \cite[Lem 9.6]{hashorva2021shiftinvariant} the local RF $\Theta^*$ defined by $Z^*$  is stochastically continuous and hence $Q^*$ constructed by $\Theta^*$ is stochastically continuous.  
\QED

\proofkorr{korEVG} \underline{Proof of \Cref{enmP1}}:  Under condition \eqref{BR2} we have that $\kk$ is purely dissipative and hence the claim follows since $Z$ is almost surely positive.\\
\underline{Proof of \Cref{enmP2}}: For our construction both $W_i$ and $W_i^{(\alpha)}$ drift to $\IF$, hence in view of \Cref{kabA} and \Cref{conv} $\kk$ has a random-shift representation. 
\QED 

\proofprop{MAPROX} 
Since $K$ is bounded, given a fixed  $t_0 \in \R^l$,  there exists a positive integer $k$
 such that  $K- t_0 \in [0,k]^l$.  
 Hence by the definition, the shift-invariance of the measure $\lambda(\cdot)$, the $\alpha$-homogeneity of $\kappa$, the assumption that $\kappa(0)=0$  and the representations for $\vk Z_N$ and $\vk Z^{(m)}_N$ imply for all $n>0,m>0$  
\bqny{
\lefteqn{\frac{1}{n^l} \E*{ \sup_{ t\in n K \cap \TT} \abs{\kappa(B^{-t}\vk Z_N) - \kappa(B^{-t}\vk Z^{(m)}_N)} }}\\
&=&
\frac{1}{n^l} \intT  \E*{ \sup_{ t\in n K \cap \TT } \abs{\kappa(B^{-s-t} \vk Q)  - \kappa(B^{-s-t}\vk Q^{(m)})}}  \lambda(ds) \\
&=&
\frac{1}{n^l} \intT  \E*{ \sup_{ t\in n (K- t_0) \cap \TT}\abs{ \kappa(B^{-s-t} \vk Q ) - \kappa(B^{-s-t}\vk Q^{(m)}) }}  \lambda(ds) \\
&\le &	  \frac{1}{n^l} \intT  \E*{ \sup_{ t\in [0,nk]^l } \kappa(B^{-s-t} \vk Q )  
\ind{  \norm{t-s} >m} } \lambda(ds) \\	   
&\le &	  \intT  \E*{ \sup_{ t\in [0,k]^l } \kappa(B^{-s-t} \vk Q )  
	\ind{  \norm{t-s} >m} } \lambda(ds) \\	
&	 \to 0, &  \quad m\to \IF,
}
where the finiteness of the integral in the second last line above follows from \eqref{dense}, which in particular implies that 
$\kappa(B^{-t} Q)$ is a.s.\ finite for all $t\in \TT$. Consequently, since further $\kappa(0)=0$ we have that a.s.\ 
\bqny{ \kappa(B^{-s-t} Q)- \kappa(B^{-s-t} Q\ind{\norm{ t-s} \le m})&=&
\kappa(B^{-s-t} Q)- \kappa(B^{-s-t}Q ) \ind{\norm{ t-s} \le m}\\
&=& \kappa(B^{-s-t} Q)\ind{\norm{ t-s} > m}
}
establishing the proof. 
	\QED

\proofprop{BRprop} For $\kappa$ as in \Cref{examp:kappa}, \Cref{examp:kappa3} the fact that the matrix pseudo-cross variogram function $\gamma$ defines uniquely 
the local RF $\Theta$ and thus also $\mathcal{C}_{\kappa}$ is shown in \cite[Example 4.2]{hashorva2021shiftinvariant}. In fact, this claim follows also 
\cite[Lem 4.2]{KumeE}. Consequently, from the statement of \Cref{thetachanges} we have that the law of the local RF $\Theta_\star$ or $\mathcal{C}_{\kappa_\star}$ depends only on the law of $\Theta$. Since $Q$ can be defined by $\Theta_\star$, it follows that also the law of $Q$ is determined only by $\kappa_\star$ and $\gamma$  establishing the claim. 
\QED 

\section{Technical Results}\label{appn}
\BEL 
Let  $\vk U\in \clasP$ and $N$ the $\TTT$-valued rv be defined on the complete probability space $\Pspace$. If $N$ is independent of $U$,  then 
$B^N \vk U\in \clasP$. If further $ {U}$ is stochastically continuous, then so is $ {B^N U}$.
\label{molla}
\EEL  
\prooflem{molla} We show first that  $B^N \vk U(t)= \vk U(t- N),t\in \TTT$ is a well-defined RF. If $N$ is a discrete rv taking only finite values $t_i, 1 \le i\le n$, then for any $c\inr$ we have 
$\{\vk U(t- N) < c \}$ is an event since $\{\vk U(t- N) < c , N=t_i\}$ is an event for all $1 \le i \le n$. For a general $N$ we can approximate it a.s.\ by discrete rvs $N_n, n\inn$. Hence 
$   B^N \vk U(t), t\in \TTT$ is the a.s.\ limit of $B^{N_n} \vk U(t), t\in \TTT$ as $n\to \IF$.
The a.s.\ limit of jointly measurable and separable  RFs  with separant $\diad$ is clearly measurable and separable with separant $\diad$.  Since $N$ is independent of $\vk U$  by the dominated convergence theorem, we have that if $ { U}$ is stochastically continuous, then $ {B^N U}$ is also stochastically continuous establishing the claim. 
\QED

\def\EJ{0}
\def\AAA{\mathcal{A}}
\def\AAAT{\mathcal{\chi}_1(\vk \Theta)}
\def\AAAY{\mathcal{\chi}_2(\vk Y)}
\def\AAAZ{\mathcal{\chi}_3(\vk Z)}
\def\EIE{\mathcal{V}}
\def\EI{\mathcal{L}}

Recall that in our notation $\E{ A;B}$ stands for  $\E{A \ind{B}}$.
Hereafter $\EIE\subset \EI$ are two additive subgroups  of $\TT$ with $\EI$ having countably infinite number of elements and set below 
$$ 
\mathcal{J}^\star_k[f]= \sum_{i \in \EIE}   \ind{\mathcal{J}_k(B^{-i}f)=0 } , \quad  k=1,2,3.
$$
\BEL 
Let $F \in \Hh_0, \Gamma \in \TwH$.  If   $\mj_1$ is a  positive shift-involution, $\mj_2$ is anchoring and $\mj_3$ is a shift-involution, respectively, then  
\bqn{
\EHh*{ F( \vk \Theta  ) \JFA [\Theta]   } 	&=&	\label{e11AB} 
\EHh*{ \sum_{i\in \TTd} \kappa(B^i\vk \Theta) F(B^{i} \vk \Theta)  ; \mj_2(\vk \Theta) =0 }, \\
\label{e11A}  
\EHh*{  \Gamma( \vk Y ) \JFB[Y] }  &=&  
\EHh*{ \sum_{i \in \TTd}    \ind{ \kappa(B^i  \vk Y)>1 } \Gamma(B^{i} \vk Y) ;  \mj_2(  \vk Y  ) =0 }  ,  
\label{gjelle}\\
\E*{ \kappa(\vk Z)  F( \vk Z )\JFC [Z]   } &=&
\E*{ \sum_{i\in \TTd}  \kappa( B^i\vk Z) F(B^{i} \vk Z) ;  \mj_3 (\vk Z) =0  } \label{gjelle2}
}
and  if $\EIE$ is also a lattice on $\TT$ with infinite number of elements we have further
\bqn{
\quad \quad  	\quad 	\pkh{\mathcal{S}_\EIE( \vk Y) = \IF, \cED{
	\JFA[\Theta]=1 }}=	\pkh{\mathfrak{B}_{\EIE, \tau}(\vk Y) = \IF, \cED{\JFB[\vk Y]=1 }}
=\pkh{\mathcal{S}_\EIE(\vk Z) = \IF,  \cED{\JFC [Z]=1 }}=0
\quad \label{0eing}
}
for all $\tau \inr$.
\label{0LZR}
\EEL

\BRM 
\begin{enumerate}[(i)]
\item 
If \eqref{lere} holds and further  
$$
\TTd= \EI=\TT=\Z^l, \quad \pkh{ \JFB[ Y]=1}=1, $$
 then \eqref{gjelle} reduces to  \cite[Prop 3.6]{BP} and \cite[Prop 3.2, Eq.\ (17)]{Planic}.
\item Let 
$\overline{\vk \Theta}(t)=\vk \Theta(t)/\sup_{t \in \EIE} \kappa^{ 1/\alpha}(B^{-t}\vk \Theta)$ 
and recall that $\E{\kappa( \vk Z)}=1$. Taking  
$\mj_k, k\le 3$  as in \nelem{0LZR} 
satisfying  \eqref{peres2}  and \eqref{peres1}, for all   $F\in \Hh_\alpha  $ shift-invariant with respect to $\EI$ 
\bqn{ \label{0e11A} 
\EHh*{  \frac{F(\vk \Theta  )}{\STq} } =
\EHh*{  F(\vk \Theta)  ; \mj_1(\vk \Theta) =0  }
=	
\EHh*{  F(\overline{\vk \Theta})  ; \mj_2(\vk Y) =0  } =		
\E*{  F(\vk Z) ;  \mj_3 (\vk Z)=0  } 
} and when $\pk{\ST< \IF}=1$ and hence by \eqref{arm10} also $\pk{\mathcal{J}_2^*[Y]=1}=1$ 
for all $\Gamma=\TwH$   shift-invariant with respect to $\EI$ and all $\tau \inr ,b\ge 1  $   
\bqn{\quad 
		b^\alpha
	\EHh*{\frac{  \kappa^\tau(Y)  \ind{\MPY> b  } \Gamma(Y)}{ [\MPY]^\alpha \mathfrak{B}_{\EI, \tau}( \vk Y)  } } 
&=&\EHh*{ \frac{ \kappa^\tau (Y)\Gamma( \vk Y )}{ [\MPY]^\alpha\mathfrak{B}_{\EI, \tau}( \vk Y) }   }  = 
\EHh*{  \Gamma(\overline{\Theta}) ;  \mj_2(  \vk Y  ) =0 
}\label{0gjelle}
\\
&=& b^\alpha 	\EHh*{  \Gamma(\overline{\Theta})\ind{\MPY> b  } ;  \mj_2(  \vk Y  ) =0 
},\label{0gjelle:2}
}
where the first equality follows from \eqref{deshtimos2} and the last is consequence of the first and the third; we interpret $\IF/\IF$ and $0/0$ as 0.  Note in passing that \eqref{0eing}  and both 
\eqref{0e11A},  \eqref{0gjelle} imply for $\TT=\Z^l$ the claims of \cite[Thm 5.5.3]{kulik:soulier:2020}. 
The first two identities in \eqref{0e11A} as well as \eqref{0gjelle}  are stated  in \cite{kulik:soulier:2020} for $\tau=0, \TTd=\EI=\TTT= \Z^d$.		 
\end{enumerate}  
\ERM 

\prooflem{0LZR} If $\mj_1$ is a positive shift-involution, then for all $F \in \Hh_0$ 
\bqny{
	\EHh*{ \sum_{i\in \EIE}  \kappa(  B^i \vk  \Theta) F( B^i \vk \Theta);    \mj_1(\vk \Theta )  =0  }
&=& \sum_{i\in \EIE} 
\EHh{    F(  \vk \Theta) \mathbb{I}(   
\mj_1(B^{-i}\vk \Theta) =0 ) \mathbb{I} (\kappa( B^{-i}\vk \Theta) \not  =0  ) }\\
&=& \sum_{i\in \EIE} 
\EHh{    F(  \vk \Theta) \mathbb{I}(   \mj_1(B^{-i}\vk \Theta)  =  0) }\\
&=& 
\EHh{ F(  \vk \Theta )  \JFA[\Theta]},
}
where the first equality follows by \eqref{eqDo20}, whereas the second one follows by  \Cref{A3}, 
hence \eqref{e11AB} follows. 

Applying  \eqref{tYY}, 
 for all $ \Gamma \in \TwH$  and $\mj_2$  anchoring (hence $\mj_2( B^{-i}\vk Y  ) =0$ implies $\kappa(B^{-i}\vk Y)> 1$ used below to derive the second equality) 
\bqn{ 
	\EHh*{  \sum_{i \in \EIE}    \ind{\kappa( B^i \vk Y )> 1} \Gamma(B^i \vk Y);  \mj_2(  \vk Y   ) =0 } 
&=&   \sum_{i \in \TTd}  
\EHh*{     \ind{\kappa(B^{-i}\vk Y)> 1}   \Gamma(\vk Y); \mj_2( B^{-i}\vk Y  ) =0 } \notag \\ 
&=&   \sum_{i \in \TTd}  
\EHh*{     \Gamma(\vk Y); \mj_2(B^{-i}  \vk Y  ) =0 }\notag\\
&=&     \EHh*{   \Gamma( \vk Y)  \sum_{i \in \TTd}  \ind{ \mj_2( B^{-1}\vk Y ) =  0}     } \notag\\
&=&     \EHh*{   \Gamma( \vk Y) \JFB[Y] }
\label{morges2}
}
and thus \eqref{e11A} follows.  Set  next 
$$\Gamma(\vk Y) = \ind{ \mathfrak{B}_{\EIE ,\tau}(\vk Y) =\IF}=
 \ind{ \mathfrak{B}_{\EIE ,\tau}(\vk Y) =\IF}\kappa^\tau(Y)/R^\tau $$
 a.s.,  where $Y= R \Theta$ (recall $\kappa(Y)=R$ a.s.). 
Since $\EIE$ is a  subgroup of the additive group $\TT$, then from \eqref{morges2}
\bqny{ 
	\EHh*{  \sum_{i \in \EIE}    \ind{\kappa( B^i \vk Y )> 1} \Gamma(B^i \vk Y);  \mj_2(  \vk Y   ) =0 } 
	= \EHh*{  \mathfrak{B}_{\EIE ,\tau}(\vk Y)  \ind{ \mathfrak{B}_{\EIE ,\tau}(\vk Y) =\IF} /R^\tau ;  \mj_2(  \vk Y   ) =0 } 
	&=&     \EHh*{   \Gamma( \vk Y) \JFB[Y] }
}

 Borrowing the idea of \cite{kulik:soulier:2020}, since  $\EIE$ has infinite  number of elements, from the above and \Cref{A1},  we conclude that
\bqn{\label{morges3}
\pkh*{  \mathfrak{B}_{\EIE,\tau}(\vk Y) =\IF,  \JFB[\vk Y]=1 }=0.
}
Next, taking  $\mathcal{J}_3$ to be  a  shift-involution and thus  it is 0-homogeneous, utilising  the shift-invariance of $\kk$,  for all $F \in \Hh_0 $  the Tonelli Theorem  implies 
\bqny{	 
\E*{    \sum_{i\in \EI}  \kappa(  B^i \vk Z) F( B^i \vk Z) ;    \mj_3(\vk Z ) =  \EJ } 	
&=&   \sum_{i\in \EI}  \E{   \kappa(Z)   F(   \vk Z) \mathbb{I}(   
\mj_3 (B^{-i}\vk Z)  = \EJ ) }\\
&=& \E{   \kappa(    Z)  F(    \vk Z)  \mj_3^*[Z]  }\\
&=&  \E{   \kappa(   \vk  Z)  F(  \vk Z/ \kappa^{1/\alpha}(   \vk  Z))\mj_3^*[Z/\kappa^{1/\alpha}(Z)]    }\\
&=& \E{\kappa( \vk Z) }\EHh{ F(  \vk \Theta ) \mj_3^*[\Theta ]   },
}
which proves  \eqref{gjelle2}. Since  $\EIE$  has  countably infinite number of elements, then   $F(\Theta)= \ind{\mathcal{S}_{\EIE} (\Theta) =\IF}$  is  shift-invariant with respect to $\EIE$. From the above   and \Cref{A1}
$$\E*{ \mathcal{S}_{\EIE} (\vk Z)F(\vk Z) ; \mj_3(\vk Z) = 0    }=
\EHh*{ F(\vk \Theta) \JFC[\Theta]  }= 
\pkh*{ \mathcal{S}_{\EIE} (\Theta ) = \IF,  \JFC[\Theta]=1   }.
$$
Borrowing the idea of the proof of \cite[Thm 5.5.3]{kulik:soulier:2020}, by the above choice of $F$ 
$$
\E*{ \mathcal{S}_{\EIE} (\vk Z)F(\vk Z) ; \mj_3(\vk Z) = \EJ    }\in \{ 0, \IF\}$$
and thus applying further \cite[Lem 9.7]{hashorva2021shiftinvariant} we obtain   
$$
\pk*{ \SZq = \IF,  \JFC[Z]=1 }=
\pkh*{ \STq = \IF, \JFC[ \Theta] =1  }=0,
$$
which together with \eqref{morges3} establishes  \eqref{0eing} and thus the proof is complete.
\QED

\def\EJ{\mathcal{Q}}
\def\MME{ \mathcal{M}_\EI(\vk Y)} 
\def\MMME{ \mathcal{M}_\EI(\vk Y)} 
\def\cEp#1{#1}

{\bf Acknowledgment}: 
I would like to thank Dmitry Zaporozhets, Dmitry Korshunov, Krzys D\c{e}bicki, Georgiy Shevshenko, Ilya Molchanov and Philippe Soulier for  several helpful discussions. 
  A special great thanks goes to the referee and the handling Editor  for their time, critical comments and suggestions which improved the manuscript significantly. Partial support by SNSF Grant 200021-196888 is kindly acknowledged.

{\bf Data Availability}. 	 Data sharing not applicable to this article as no datasets were generated or analysed during	the current study.

{\bf Conflict of interest}.  The author declares that he has no conflicts of interest to this work.

\bibliographystyle{ieeetr}
\bibliography{JTPresubmission.bib}

\end{document}